\documentclass[12pt]{amsart}
\usepackage{epsf} 
\usepackage{amscd}
\usepackage{graphics}
\usepackage{graphpap}
\usepackage{amsmath}
\usepackage{amssymb}
\usepackage{eucal}
\usepackage{amsfonts}
\usepackage{amsmath}
\usepackage[Lenny]{fncychap}
\usepackage{enumerate}
\usepackage{amssymb}
\usepackage[english]{babel}
\usepackage[all]{xy}
\usepackage{fancybox}
\usepackage{latexsym,mathrsfs,longtable,lscape,dsfont,yfonts}
\usepackage{srcltx}
\usepackage{soul}

\usepackage[normalem]{ulem}

\def\cc {{\mathfrak c}} 
\def\AA {{\mathbb A}}

\def\NN {{\mathbb N}}

\def\PP {{\mathbb P}}
\def\QQ {{\mathbb Q}}
\def\RR {{\mathbb R}}
\def\TT {{\mathbb T}}

\def\ZZ {{\mathbb Z}}

\def\sB {{\mathcal B}}
\def\sC {{\mathcal C}}

\def\sE {{\mathcal E}}
\def\sF {{\mathcal F}}

\def\sK {{\mathcal K}}
\def\sL {{\mathcal L}}
\def\sM {{\mathcal M}}

\def\sP {{\mathcal P}}

\def\sS {{\mathcal S}}
\def\sT {{\mathcal T}}
\def\sU {{\mathcal U}}

\def\bp {{\mathbf p}}
\def\bq {{\mathbf q}}

\def\bu {{\mathbf u}}

\def\cf {\mathrm{cf}}

\def\sup {\mathrm{sup}\:} 
\def\max {\mathrm{max}\:} 

\def\implies {\Longrightarrow}
\def\empty {\emptyset}

\def\and {\wedge}
\def\qed {{\blacksquare}}
\def\onto {\,\rule[.04in]{.15in}{.01in}\kern-4pt\raise.3ex\hbox{$\scriptscriptstyle>\!>\,$}}
\def\to {\longrightarrow}

\def\st {such that}
\def\tg {topological group}
\def\sii {if and only if}

\def\wrt {with respect to}
\def\tb {totally bounded}

\def\nhd {{neighborhood}}

\def\ie {{\em i.e.,}}

\def\< {{\langle}}
\def\> {{\rangle}}

\def\eps {{\epsilon}}

\def\gs {{\sigma}}

\def\^ {{\widehat{ \:}}}
\def\wZ {{\widehat{\ZZ}}}
\def\wG {{\widehat{G}}}
\def\wN {{\widehat{N}}}

\def\wH {{\widehat{H}}}

\def\pf {{\em \noindent Demostraci\'on:}}
\def\epf {~\hfill $\qed$}

\def\pf {{\em \noindent Proof:}} 
\def\l( {{\left)}} 
\def\r( {{\right)}} 
\def\l[ {{\left[}} 
\def\r] {{\right]}} 
\def\l{ {{\left{}} 
\def\r} {{\right}}} 

\newcommand{\mkp}{\medskip}
\newcommand{\bkp}{\bigskip}

\newfont{\cyr}{wncyr8}
\newfont{\cyb}{wncyr8}

\usepackage{pdfpages}
\usepackage{amsxtra}
\newtheorem{thm}{Theorem}[section]
\newcommand{\bthm}{\begin{thm}}
\newcommand{\ethm}{\end{thm}}

\newtheorem{prop}[thm]{Proposition}
\newcommand{\bprp}{\begin{prop}}
\newcommand{\eprp}{\end{prop}}

\newtheorem{fact}[thm]{Fact}
\newcommand{\bfct}{\begin{fact}}
\newcommand{\efct}{\end{fact}}

\newtheorem{prob}[thm]{Problem}
\newcommand{\bprb}{\begin{prob}}
\newcommand{\eprb}{\end{prob}}

\newtheorem{quest}[thm]{Question}
\newcommand{\bqtn}{\begin{quest}}
\newcommand{\eqtn}{\end{quest}}

\newtheorem{lem}[thm]{Lemma}
\newcommand{\blem}{\begin{lem}}
\newcommand{\elem}{\end{lem}}

\newtheorem{claim}[thm]{Claim}
\newcommand{\bclm}{\begin{claim}}
\newcommand{\eclm}{\end{claim}}

\newtheorem{cor}[thm]{Corollary}
\newcommand{\bcor}{\begin{cor}}
\newcommand{\ecor}{\end{cor}}

\newtheorem{conj}[thm]{Conjecture}
\newcommand{\bcnj}{\begin{conj}}
\newcommand{\ecnj}{\end{conj}}

\theoremstyle{definition}
\newtheorem{defn}[thm]{Definition}
\newcommand{\bdfn}{\begin{defn}}
\newcommand{\edfn}{\end{defn}}

\newtheorem{spec}[thm]{Specializing}
\newcommand{\bspc}{\begin{spec}}
\newcommand{\espc}{\end{spec}}

\theoremstyle{remark}
\newtheorem{rem}[thm]{Remark}
\newcommand{\brem}{\begin{rem}}
\newcommand{\erem}{\end{rem}}

\newtheorem{cnv}[thm]{Convention}
\newcommand{\bcnv}{\begin{cnv}}
\newcommand{\ecnv}{\end{cnv}}

\newtheorem{exam}[thm]{Example}
\newcommand{\bexm}{\begin{exam}}
\newcommand{\eexm}{\end{exam}}

\newtheorem{exercise}[thm]{Exercise}
\newcommand{\bexr}{\begin{exercise}}
\newcommand{\eexr}{\end{exercise}}

\newtheorem{thmy}{\textbf{Theorem}}
\newenvironment{thmx}{\stepcounter{thm}\begin{thmy}}{\end{thmy}}

\begin{document}

\title[On closed subgroups]{On closed subgroups of precompact groups}
\date{\today}
\author{Salvador Hern\'andez}
\address{Departmento de Matem\'aticas\\
Universitat Jaume I, 8029-AP\\
Castell\'on, Spain}
\email{hernande@mat.uji.es}
\author{Dieter Remus}
\address{Institut f\"ur Mathematik\\
Universit\"at Paderborn, Warburger Str. 100 \\
D-33095 Paderborn, Germany}
\email{remus@math.uni-paderborn.de}
\author{F. Javier Trigos-Arrieta}
\address{Department of Mathematics\\
California State University, Bakersfield\\
Bakersfield, CA 93311}
\email{jtrigos@csub.edu}
\thanks{\emph{Corresponding author}: F. Javier Trigos-Arrieta}


\begin{abstract}
It is a Theorem of W.~ W. Comfort and K.~ A. Ross that if $G$ is a subgroup of a compact Abelian group, and $S$ denotes those continuous homomorphisms from $G$ to the one-dimensional torus, then the topology on $G$ is the initial topology given by $S$. {Assume that $H$ is a subgroup of $G$. We study how} the choice of $S$ affects the topological placement and properties of $H$ in $G$. Among other results, we have {made significant} progress toward the solution of the following specific questions:
How many totally bounded group topologies does $G$ admit such that $H$ is a closed (dense) subgroup?
If $C_S$ denotes the poset of all subgroups of $G$ that are $S$-closed, ordered by inclusion, does $C_S$ have a greatest (resp. smallest) element?
We say that a totally bounded (topological, resp.) group is an \textit{SC-group} (\textit{topologically simple}, resp.) if all its subgroups are closed
(if $G$ and $\{e\}$ are its only possible closed normal subgroups, resp.) {In addition, we investigate the following questions.} How many SC-(topologically simple totally bounded, resp.) group topologies
does an arbitrary Abelian group $G$ admit?
\end{abstract}

\dedicatory{Dedicated to Mar\'{\i}a Jes\'us Chasco on the occasion of her 65th birthday }
\maketitle

\section{Introduction}
Let $(G,\tau)$ be an Abelian Hausdorff topological group.
In their 1964 seminal paper \cite{ComfortRoss1964}, Comfort and Ross assigned a Hausdorff group topology $\tau_S$ to each point-separating subgroup $S$ of the character group $\wG$ of $G$, consisting of all group-homomorphisms from $G$ into the unit circle $\TT$ as follows:
$\tau_S$ is the weakest topology on $G$ that makes the elements of $S$ continuous. (By \cite[Theorem 1.9 ]{ComfortRoss1964} a subgroup of the compact group $\wG$ is point-separating \sii \ it is dense.)
Let $S_\tau:=\{\phi \in \wG: \phi \mbox{ is $\tau$-continuous}\}$. Then $S_\tau$ is a subgroup of $\wG$. By \cite[Theorem 1.2 and Theorem 1.3]{ComfortRoss1964} the following holds: (i) $(G,\tau)$ is totally bounded \sii \ $\tau=\tau_{S_\tau}$; (ii) if $S$ is a point-separating subgroup
of $\wG$, then $S_{\tau_S}=S$. (We call this assertion {\em the Comfort-Ross Theorem.})
Let $\sB(G)$ be the set of all totally bounded group topologies on $G$. By the Comfort-Ross Theorem there is an order-preserving bijection
from $\sB(G)$ onto the set of all point-separating subgroups of $\wG$.

In this paper we consider $\tau_S$ for arbitrary subgroups $S$ of $\wG$
and $S_\tau$ if $\tau$ is not necessarily Hausdorff. For an Abelian group $G$ let $\sP\sK(G)$ be the lattice of all precompact group topologies on $G$, and let $\Sigma(\wG) $ be the lattice of all subgroups of $\wG$. By using the Comfort-Ross Theorem, Remus in \cite{Remus1983} showed that
$f: \sP\sK(G)\to\Sigma(\wG) $ is a lattice-isomorphism, where $f(\tau)=S_\tau$ and $f^{-1}(\mu)=\tau_\mu$ hold. (For a generalization to arbitrary groups see \cite[(3.7)]{Remus1983} and \cite{Remus1985}.)

Fast forward to 1983 and 1985 when Remus in \cite{Remus1983} (see also \cite{Remus1986}) and  Berhanu, Comfort and Reid in \cite{BerhanuComfortReid1985} proved that if $G$ is infinite, then $\wG$ has $2^{2^{|G|}}$-many dense subgroups. It follows that an infinite group accepts $2^{2^{|G|}}$-many \tb \ (precompact and Hausdorff) group topologies. Concerning the number of totally bounded (in particular pseudocompact) group topologies on non-necessarily Abelian groups we refer to \cite{Remus1991}, \cite{ComfortRemus1994a}, \cite{RemusNovember1995}, \cite{ComfortRemus1996}, \cite{ComfortRemus1997} and \cite{ComfortRemus2016}.

Many examples regarding the interplay between precompact or \tb \ group topologies between $G$ and subgroups of $\wG$ have been studied by several colleagues. The following list is by no means complete or exhaustive:

\begin{enumerate}
\item
The authors of \cite{comfsound82} showed that if
$\tau_1 \subsetneq \tau_2$ are (Hausdorff) group topologies on $G$ \st \ $(G,\tau_1)$ is compact and $(G,\tau_2)$ is pseudocompact,
then there are $\tau_2$-closed subgroups of $G$ that are not $\tau_1$-closed. They also proved that if $(G,\tau_1)$ is a totally disconnected Abelian compact group of uncountable weight, then there is a pseudocompact group topology $\tau_2$ \st \ $\tau_1 \subsetneq \tau_2$. Eventually,  the authors in \cite{ComfortRobertson1987} generalize this result by removing the requirement of $(G,\tau_1)$ being totally disconnected.
	
\item The authors of \cite{jana_sound82} focus on an infinite compact (Hausdorff) totally disconnected Abelian group $(G,\tau)$ and try to
obtain finer totally bounded group topologies $\tau'$ such that every $\tau'$-closed subgroup is $\tau$-closed.
	
\item In \cite{soundaiii}, the author focuses on \tb \ \tg s in which every subgroup is closed.

\item In \cite[Proposition 3.4]{HernandezMacario2003}, the following is proved: Let $G$ be a \tb \ Abelian group with character group $S$. If $L$ is a subgroup of $S$, let $L_G$ denote $L$ equipped with the weakest topology that makes the elements of $G$ (acting on $L$) continuous. Then, $G$ is pseudocompact if and only if the topology inherited by each countable subgroups of $S_G$ is its cor\-res\-pon\-ding largest \tb \ group topology.
\item In \cite{CDT}, the following is proved: Let $G$ be a precompact, bounded torsion Abelian group with character group $S$. If $G$ is Baire (resp., pseudocompact), then all compact (resp., countably compact) subsets of $S_G$ \ are finite. Also, $G$ is pseudocompact if and only if all countable subgroups of $S_G$ are closed.
\end{enumerate}

In this paper we further investigate the topological properties of precompact and totally bounded abelian groups via its dual group.
More precisely, if $H$ is a subgroup of a precompact Abelian group $G$, what are the topologies of the form $\tau_S$, with $S$ a subgroup of $\wG$ \st \

\begin{enumerate}
	\item $H$ is $\tau_S$-closed?
	\item $H$ is $\tau_S$-dense?
	\item How many subgroups $S$ of $\wG$ are there \st \ each of the above happens?
	\item How many subgroups of $\wG$ produce the same closed (dense) subgroups in $G$?
\end{enumerate}

Similarly, we want to know those subgroups of $\wG$ producing the same closed (dense) subgroups in $G$.\mkp

We now formulate our main results.\mkp

\begin{thmx}\label{Thm_A}
\it{If $G$ is an infinite Abelian group and $\{H_i: i\in I\}$ is a family of subgroups of $G$ such that $|I|<2^{|G|}$, then $G$ admits exactly $2^{2^{|G|}}$-many totally bounded group topologies $\mu$ such that $H_i$ is closed in $(G,\mu)$ for all $i\in I$}.
\end{thmx}
\mkp

{By Proposition \ref{khan1} this result is not true if $|I| = 2^{|G|}$.}

{Let $G$ be an Abelian group and let $S$ be a subgroup of $\wG$. We denote by $C_S$ the poset of all subgroups of $G$ that are $S$-closed.
Hence $C_S$ is the poset of all subgroups of $G$ which are closed in the precompact group $(G,\tau_S)$.}

{It is natural to ask whether there exists a greatest precompact group topology $\tau_{MS}$ on $G$ with $C_{MS}=C_S$. The next result gives a positive
answer, since $\sP\sK(G)$ and $\Sigma(\wG)$ are isomorphic as lattices.}

\begin{thmx}\label{Thm_B}
\it{Let $G$ be an Abelian group and let $S$ be a subgroup of $\wG$. Then there exists a greatest subgroup $M{S}$ containing $S$ and \st \ $$C_{MS}~=~C_S.$$}
\end{thmx}
\mkp

We remember that a totally bounded Abelian group is an \textit{SC group} if all its subgroups are closed (see Remark \ref{newrem 1}.
Although this seems to be a very restrictive
property, the following result shows that if $G$ is an Abelian group that is not of bounded order, then the number of SC-group topologies is huge.
\bkp

\begin{thmx}\label{Thm_C}
\it{The following statements hold:

\begin{enumerate}
  \item[(a)] Let $G$ be an Abelian group which is not of bounded order. Then $G$
admits at least $2^\cc$-many SC-group topologies.
  \item[(b)] Every countable Abelian group which is not of bounded order admits
exactly $2^\cc$-many SC-group topologies.
\end{enumerate}}
\end{thmx}
\bkp

We remember {that a topological} group $(G,\tau)$ is \emph{topologically simple} if $G$ and $\{e\}$ are its only possible closed normal subgroups.
The anti-discrete topology on a space $X$ is defined as $\{\empty, X\}$. Hence, if $G$ is a topologically simple  \tg,
then $G$ is Hausdorff unless $G$ carries the anti-discrete topology (\cite[(5.4)]{hewitt1963abstract}).

By Corollary \ref{ts2} every infinite Abelian totally bounded group which is topologically simple is algebraically a subgroup of $\RR$.
The next results clarifies the question of the existence of topologically simple for subgroups of the real line.
\bkp

\begin{thmx}\label{Thm_D}
\it{Let $G$ be a non-trivial subgroup of $\RR$. Then $G$ admits exactly $2^\cc$-many \tb \ group topologies $\tau$ such that $(G,\tau)$ is topologically simple. The topologies $\tau$ can be chosen such that $w(G,\tau)=\cc$.}
\end{thmx}

To prove the foregoing, we develop some technical results we believe are interesting on its own. For example, we characterize the subgroups
of the compact character group $\wG$ of an Abelian group $G$, producing the same closed subgroups in $G$. This {is Theorem} \ref{5.3}.
\mkp

{This paper is organized as follows: Sections 1, 2 and 3 consist of the introduction, notation and results about dense subgroups of precompact groups.
In Section 4, we are concerned with the number of totally bounded
group topologies on an Abelian group $G$ such {that  each member of a previously fixed family of subgroups of $G$ is closed
on all those topologies.} The goal of Section 5 is {somewhat different as we deal with the \emph{poset} of all subgroups $S$ of the \emph{big} dual group $\wG$
having a previously fixed subgroup $H$ of $G$  closed} in the weak topology $\tau_S$. In Section 6, we characterize and then calculate
the number of totally bounded group topologies on $G$ such that all subgroups of $G$ are closed.
In Section 7, we go in the opposite direction {as we investigate the existence of topologies without producing closed non-trivial subgroups.}
In Section 8, we apply previous results to one specific example: The integers. 
{In Section 9 we conclude with final remarks.}}

\section{Notation}

By $\ZZ$ we denote the group of integers and by $Z(n)$ the cyclic group of order $n$. Sometimes we look at $\ZZ$ as a topological, discrete group, and some others as the underlying group of it as a \tg \ but without any topology. Set $\omega:=\{0,1,2,...\}$ and
$\NN:=\{1,2,...\}$. $\PP$ denotes the set of prime numbers. {Our model for $\TT$ is the group $([0,1),+$~mod~$1$)} with the topology inherited from $\RR$, when we need to see it as a \tg . Let $G$ be an Abelian group. If $A \subseteq G$, the subgroup generated by $A$, namely $\< A \> $, is the smallest subgroup containing $A$; if $A$ is the singleton $\{a\}$, we write just $\< a \> $. When writing $H \leq G$, we signify that $H$ is a subgroup of $G$. The symbols $r(G), r_0(G)$ and $r_p(G)$ stand respectively for the {\em rank, torsion-free rank,} and {\em $p$-rank} of the group $G$ \cite[\S 16]{fuchsi}. $tG$ denotes the torsion subgroup of $G$.
{If $L$ and $M$ are groups which are algebraically isomorphic, then we write $L\cong M$. If the topological groups $G$ and $H$ are topologically  isomorphic, we write $G\simeq H$.} For an Abelian group $G$ we will denote by $\wG$ the set of all  homomorphisms $\phi:G \to \TT$, which we will also refer to as {\em the characters of} $G$. $\wG$ becomes a group by defining $(\phi_1\phi_2)(g):=\phi_1(g)+\phi_2(g) \in \TT$ whenever $g\in G$, and equipped with the finite-open topology $\gs(\wG,G)$, $\wG$ becomes a compact \tg . We know $\wZ\simeq \TT$.
When $(G,\tau)$ is a topological Abelian group, then $(G,\tau) \^ $, or simply $G\,\^ \ $ if there is no room for confusion, denotes the subgroup of $\wG$ consisting of the $\tau$-continuous elements.
If $S$ is a subgroup of $\wG$, denote by $G_S$ the \tg \ obtained by equipping $G$ with the weakest topology $\tau_S$ that makes the elements of $S$ continuous. It follows that $G_S$ is Hausdorff \sii \ $S$ separates the elements of $G$ (\ie \ $S$ is {\em point-separating}.)
We have pointed out above that the Comfort-Ross Theorem establishes an order-preserving bijection between the set of all to\-tal\-ly bounded group topologies on $G$
and the set of all point-separating sub\-groups of $\wG$. Therefore, taking $S=\wG$, yields the finest precompact group topology on $G$. This topology
is called the \emph{Bohr topology} of $G$ and it is de\-sig\-na\-ted by $\tau_b(G)$ here. It is a well-known fact that the Bohr topology is Haus\-dorff (see \cite{comf-hbook}).

{As usually, if $Y$ is a subspace of the topological space $X$ we let $\overline{Y}^X$ de\-no\-te the closure of $Y$ in $X$;
$wX$ and $\chi X$ stand respectively for the {\em weight} and {\em cha\-rac\-ter} of $X$ \cite[\S 3]{comf-hbook}. Also, given a subgroup $S$ of $\wG$, we {let $H_{|S}$ denote the} group $H$ equipped with the topology inherited from $G_S$, and $\AA(S,H):=\{\phi \in S:\phi[H]=\{0\}\}$. $\AA(S,H)$ is called {\em the annihilator of $H$ in $S$} and is a subgroup of $S$. Similarly, $\AA(H,S):=\{g\in H:\varphi(g)=0\:\forall \varphi \in S\}$ is called {\em the annihilator of $S$ in $H$} and is a subgroup of $H$.}

A \tg \ $G$ is {\em precompact} if whenever $U$ is an open subset of $G$, there is a finite subset $F\subset G$ \st \ $G=FU$. If in addition to be precompact, $G$ is Hausdorff, then we say that $G$ is {\em \tb .} It is a Theorem of A. Weil \cite{weil1937} that the completion of the totally bounded group $G$ is a compact group $\overline{G}$, which we call its {\em Weil completion.} (The notion \textit{compact space} is used as in
\cite[Definition 1, p.~83]{Bourbaki1966}.)

Given a \tg \ $G$, it follows that the closure $N$ of its identity is a (normal) subgroup of $G$ \cite[(5.4)]{hewitt1963abstract}, hence $G/N$ is a Hausdorff \tg \ \cite[(5.21)]{hewitt1963abstract}. We refer to the map $\phi:G\to G/N$, or simply to $G/N$, the {\em Hausdorff modification of $G$.} When the topology of $G$ stems from a subgroup $S$ of $\wG$, it follows that $N=\AA(G,S)$ if $G$ is Abelian. A topological group $G$ is precompact \sii \
its Hausdorff modification is totally bounded. By \cite[chapter III, p.~248]{Bourbaki1966} this means that the {\em Hausdorff completion} of
$G$ is compact.

Often we will have to use the following.

\blem \label{2.1} Let $G$ be an infinite Abelian group, $H$ a subgroup of $G$ and $S$ a subgroup of $\wG$. Then  $\widehat{G_S/H}$ and $\widehat{H_{|S}}$, resp., are group-isomorphic to $\AA(S,H)$ and $S/\AA(S,H)$, resp.
\elem

\pf \ We sketch the proof that $\widehat{G_S/H}$ is group-isomorphic to $\AA(S,H)$. Set $N:=\AA(G,S)$. We then have

\[\begin{CD}
H	@>>>  NH/N    @>>>     Y\\
@VVV      		 @VVV	@VVV\\
G_S 		 @>>> 	G_S/N @>>> \Gamma\\
@VVV      		 @VVV	@VVV\\
G_S/H	@>>>  G_S/NH    @>>>     \Gamma/Y\\
\end{CD}\] where the groups in the middle column are the Hausdorff modifications of the groups on the left column, $\Gamma$ is the Weil completion of $G_S/N$ and $Y$ the closure of $NH/N$ in $\Gamma$. The first vertical arrows 
are containments whereas the second vertical arrows 
are projections.
It follows that $\AA(S,H)=\AA(S,NH)$, $\Gamma$ is a compact Hausdorff group and each of the groups in the middle column is dense in the corresponding group on the right column. {To see, for example, that $\AA(S,NH/N)\cong \AA(S,NH)$, one readily sees that the map $\phi \mapsto \phi \circ \pi$ is an isomorphism, where $\pi:NH \to NH/N$ is the natural map. We then have that $\widehat{G_S/H} \cong \widehat{G_S/(NH)}\cong \widehat{\Gamma/Y}\cong  \AA(S,Y)$ (by \cite[(24.5)]{hewitt1963abstract}) $\cong \AA(S,NH/N)\cong \AA(S,NH)\cong \AA(S,H)$.}

That $\widehat{H_{|S}}$ is group-isomorphic to $S/\AA(S,H)$ can be done in a similar fashion \epf

\section{Dense Subgroups}

In the following let $G$ be an Abelian group.
Let $H$ be a subgroup of $G$. Obviously, if $H$ is dense as a \tg \ in $G_S$, then $\varphi[H]$ will be dense in  $\varphi[G] \subseteq \TT$, whenever $\varphi\in S$. We would like to prove that the latter condition is also sufficient.

\begin{thm} \label{Th1} Let $H$ be a subgroup of $G$ and $S$ a 
	subgroup of $\wG$. Then $H$ is dense in $G_S$ \sii \ $\varphi[H]$ is dense in  $\varphi[G] \subseteq \TT$, whenever $\varphi\in S$.
\end{thm}

\pf \ We must show that if $V$ is an open \nhd \ in $G_S$, there exists $h \in H \cap V$. Let $g \in V$. By definition, there exist $\phi_1,...,\phi_k \in S$ and $\eps > 0$ \st \
\[\cap_{j=1}^k \phi_j^{-1}[V_\eps(\phi_j(g))] \subseteq V,\]
where $V_\eps(\zeta):=\{t\in \TT:|\zeta-t| < \eps\}.$
Consider $\overline{S_H}$, the completion of the Hausdorff modification of $S$ equipped with the weakest topology that makes the elements of $H$ continuous, the latter viewed as characters on $\wG$. Let $\phi_1, \phi_2 \in S$. If $\phi_1(h)=\phi_2(h)$ for all $h\in H$, then $\phi_1\phi_2^{-1}(h)=0$, {which implies $\phi_1\phi_2^{-1}[H]=\{0\}$}. Since $\phi_1\phi_2^{-1}\in S$, our hypothesis implies that $\phi_1\phi_2^{-1}[G]=\{0\}$, hence we would have $\phi_1=\phi_2$. Therefore, {$\phi_1\neq \phi_2$ implies
there is $h \in H$} \st \ $\phi_1(h)\neq \phi_2(h)$, and thus it follows that $\overline{S_H}$ is compact. By \cite[(26.16)]{hewitt1963abstract}, given $\phi_1,...,\phi_k \in S_H \subseteq \overline{S_H}$, and $g\in G$, there exists $h \in H$ \st \
\[|\phi_j(h)-\phi_j(g)| < \eps, \:\:\:\:\:\:\:(j=1,2,3,..., k).\]
Hence $h \in \cap_{j=1}^k \phi_j^{-1}[V_\eps(\phi_j(g))] \subseteq V$, as required. \epf
\mkp

{Let $S$ be a subgroup of $\wG$. Obviously, $H_{|S}$ is closed in $G_S$ \sii \ whenever $H_{|S}$ is dense in $N_{|S}$, then $H_{|S}=N_{|S}$,
for all $H\leq N\leq G$. Set $T:=\{f_{|N}:f \in S\}$. Obviously $T$ is a subgroup of $\wN$; if $S$ separates the points of $G$,
then $T$ separates the points of $N$
and by Comfort-Ross' theorem \cite[Theorem 1.2]{ComfortRoss1964}, $N_{|S}=N_T$.}

\begin{lem}\label{AA}
The subgroup $H\leq G$ is closed in $G_S$ if and only if
$$\AA(\wG,H)=\overline{\AA(S,H)}^{\wG}=\overline{S\cap \AA(\wG ,H)}^{\wG}.$$
\end{lem}

\pf \ If $H$ is closed in $G_S$, then $\widehat{G_S/H}\cong \AA(S,H)$ by Lemma\ref{2.1} and $G_S/H$ is Hausdorff. Therefore $\AA(S,H)$ is dense in $\AA(\wG,H)$.

Conversely, since  $\widehat{G/H}\cong \AA(\wG,H)$, for every $g\in G \setminus H$, there is $\phi \in \AA(\wG,H)$
\st \ $\phi(g) \neq 0$ \cite[(A.7)]{hewitt1963abstract}. Therefore, if $\AA(S,H)$ is dense in $\AA(\wG,H)$, there is $\chi \in \AA(S,H)$ \st\  $\chi(g) \neq 0$.
This implies that $H$ is closed in $G_S$ by Lemma \ref{A1} .
\epf

\begin{cor}\label{A} {If $\AA(\wG,H)\subseteq S$, then $H$ is closed in $G_S$. If $H$ is of finite index in $G$, then the converse is true. }
\end{cor}

\pf \ The first assertion holds, since $\AA(\wG,H)$ is closed in $\wG$ by \cite[(23.24)(c)]{hewitt1963abstract} . Assume that $H$ is a closed subgroup of finite index in $G_S$. Then $H$ is open in $G_S$. If $\varphi \in \AA(\wG,H)$, then $\varphi_{|H}=0$ is continuous on $H$, hence continuous on $G_S$, hence $\varphi \in S$, as required.
\epf

If $H$ is not of finite index in $G$, the converse may be false:

\bexm\label{3.4}  Consider the group $G:=\oplus_\omega \< \frac{1}{4}\> $, where $ \< \frac{1}{4}\> \subset \TT$, and its subgroup $H:=\oplus_\omega  \< \frac{1}{2}\> $. Then $\wG\cong\prod_\omega Z(4)$ and $\AA(\wG,H)\cong\prod_\omega \{0,2\}$. Consider the subgroup $S:=\oplus_\omega Z(4)$ of $\wG$. Then $H$ is closed in $G_S$, yet $\AA(\wG,H)\subseteq S$ is false. [If $g \in G \setminus H$, say $g=(g_k)$, there is $n< \omega$ \st \ $g_n \not \in \< \frac{1}{2}\> $. Consider $\phi=(t_k) \in \wG$ defined as $t_k=0$ if $k\neq n$, and $t_n=2$. It follows that $\phi \in S, \phi(g)=2g_n=\frac{1}{2}$, yet $\phi[H]=\{0\}.$] \epf
\eexm

\begin{thm}\label{B} $\AA(\wG,H)\cap S=\{0\}$ \sii \ $H$ is dense in $G_S$.
\end{thm}

\pf \ {\rm($\Rightarrow$)} Deny. By Theorem \ref{Th1}, there would be $\phi \in S$ with  $\phi[H]$ not dense in $\phi[G]$. This would imply that $|\phi[H]|<\aleph_0$, hence closed in $\TT$, and there would be $g \in G \setminus H$ with $\phi(g)\not\in \phi[H]$. But then, if $\pi:\TT\to\TT/\phi[H]\simeq \TT$ denotes the canonical map, then $0\neq \pi \circ \phi \in \AA(\wG,H)\cap S$, a contradiction.

{\rm($\Leftarrow$)} Let $\phi \in \AA(\wG,H)\cap S$. By Theorem \ref{Th1} $\{0\}=\phi[H]$ would be dense in $\phi[G]$, implying that $\{0\}=\phi[G]$, hence $\phi=0$. \epf

\brem\label{maximal}
When $H$ is a {\em maximal} proper subgroup of $G$, then either $H$ is necessarily dense or closed in $G_S$.
\erem

A subgroup $E$ of an Abelian (resp. Abelian topological) group $A$ is said to be {\em essential} (resp. {\em topologically essential})
if $E \cap B \neq \{0\}$ whenever $B$ is a non-trivial (resp. and closed) subgroup of $A$ \cite[p.~84]{fuchsi}.
By \cite[Lemma 16.2]{fuchsi} an {\em independent system} $M$ is maximal \sii \ $\< M \> $ is essential.
if $M$ is a maximal independent system in an essential subgroup of $A$, then $M$ is a maximal independent system in $A$.
Note that in \cite{CD1999} and \cite{Boschi2000} topologically essential groups in our sense are called ``essential''.

The following result characterizes precompact group topologies without proper dense subgroups.
First we recall the following result \cite[(24.10)]{hewitt1963abstract}.

\begin{lem}\label{closure} { Let $S$ be a subgroup of the compact Abelian group $\wG$. Then $S$ is closed in $\wG$
if and only if $S=\AA(\wG,\AA(G,S))$. }
\end{lem}

\begin{prop}\label{no dense subgroups}  { Given $S\leq \wG$ the group $(G,\tau_S)$ contains no proper dense subgroups
if and only if $S$ is topologically essential in $\wG$. }
\end{prop}
\medskip

\pf\ Suppose that $S$ is topologically essential in $\wG$ and let $H\leq G$ be a proper subgroup of $G$. Then
$\{0\}\lneqq \AA(\wG,H)$ by \cite[(24.12)]{hewitt1963abstract}. Therefore $\AA(\wG,H)\cap S\not=\{0\}$, which implies that $H$ is not dense in $(G,\tau_S)$
by Theorem \ref{B}.

Conversely, let $L$ be a non-trivial closed subgroup of $\wG$. Then $\AA(G,L)\neq G$ by \cite[(22.17)]{hewitt1963abstract}, which
means that $\AA(G,L)$ is not dense in $(G,\tau_S)$. By Theorem \ref{B},
there is $0\not=\phi\in S\cap\AA(\wG,\AA(G,L))$ {which equals $S\cap L$ by the above result.}
This means that $S\cap L\not=\{0\}$ by Lemma \ref{closure}. \epf

\begin{cor}\label{no dense subgroups 2}  Given $S\leq \wG$ the group $(G,\tau_S)$ is a \tb\ \tg\ without proper dense
subgroups if and only if $S$ is topologically essential and dense in $\wG$.\epf
\end{cor}
\mkp

Now we are ready to give a partial answer to the problem of finding \tb \ group topologies withour proper dense subgroups.

\bthm\label{4.8} If $G$ is a torsion-free Abelian group and $S$ an essential subgroup of $\wG$, then $S$ is dense in $\wG$.
As a consequence $(G,\tau_S)$ is a \tb \ \tg \ without proper dense subgroups.

\ethm

\pf \  We know by the Comfort-Ross Theorem that $(G,\tau_S)$ is a \tb \ \tg\ if and only if $S$ is dense $\wG$. On the other hand,
in order to prove that $S$ is dense in $\wG$, it will suffice to show that if $g \in \AA(G,S)$ then $g=0$:
\noindent

Take an arbitrary element $g\in \AA(G,S)$ and define $\overline{g}:\wG \to \TT$ by $\overline{g}(\phi):=\phi(g)$.
It follows that $\overline{g}$ is a continuous character of $\wG$ \cite[(24.8)]{hewitt1963abstract}.
Since $S$ is essential in $\wG$, for every $\phi\in\wG$ there is some $m\in\NN$  \st\ $m\phi=s\in S$,
which means that $m\overline{g}(\phi)=m\phi(g)=s(g)=0$. Therefore $\overline{g}(\wG)\subseteq t\TT$. Now, since $G$
is torsion-free, it follows that $\wG$ is connected \cite[(24.25)]{hewitt1963abstract}. Thus $\overline{g}(\wG)$
is a torsion compact connected subgroup of $\TT$, which implies that $\overline{g}(\wG)=\{0\}$. In other words,
we have that $\phi(g)=0$ for all $\phi\in\wG$. This yields by \cite[(22.17)]{hewitt1963abstract} that $g=0$.

That $(G,\tau_S)$ contains no dense subgroups follows from Corollary \ref{no dense subgroups 2}.
\epf
\mkp

The converse is not true. If $S$ is the torsion subgroup of $\TT \simeq \wZ$, then $S$ is not essential in $\TT$ (see Corollary \ref{4.3} below.)

Not every $\wG$ has proper essential subgroups \cite[Corollary 16.4]{fuchsi}.

\brem\label{cex}
Theorem \ref{4.8} does not hold in general. The requirement that $G$ be torsion-free cannot be dropped: Let $m$ be an infinite cardinal, $p\in\PP$ and $n\in\NN$ with $n>1$.
For $G:=\bigoplus_m Z(p^n)$ we get $\wG=(Z(p^n))^m$. The socle of $\wG$ is $H:=\bigoplus_{2^m} Z(p)$.
Then \cite[Exercise 16.10]{fuchsi} implies that $H$ is an essential subgroup of $\wG$. The closure of $H$ is a group of order $p$.
Hence $H$ is not dense in $\wG$.
\erem

\bexm\label{4.9}  For $m\geq\omega$ and $p\in\PP$, let $G:=\bigoplus_m Z(p)$. and $S\neq \wG$. Then $G_S$ has always  non-trivial dense subgroups: For, if $\phi \not\in S$, then $G/\ker\phi\cong U$, where U is a subgroup of $t\TT$ with order $p$. Hence $\ker\phi$ is a maximal proper subgroup
of $G$. Since $\phi$ is not continuous on $G_S$, $\ker\phi$ is not closed in $G_S$. By Remark \ref{maximal} it is dense.
\eexm

On the other hand:

\bexm\label{4.9a}  If $G$ is an infinite Abelian group of bounded order, then $G_S$ has always non-trivial closed  subgroups: For, if $\phi \in S$, then $\ker\phi$ is a closed subgroup of $G_S$. Then $G/\ker\phi\cong U$, where U is a subgroup of bounded order of $t\TT$ . Hence $U$ is finite. Thus $\ker\phi$ is non-trivial.
\eexm

\bexm\label{4.9b}  If $|G|>\cc$, then $G_S$ has always non-trivial closed subgroups. For, if $\phi \in S$, then ker~$\phi$ is a closed subgroup of $G_S$ that is non-trivial since it has index at most $\cc$.
\eexm

\begin{thm} \label{Tarski-like}  Let $(G,\tau)$ be a compact Abelian group, and let $H$ be a dense subgroup of $(G,\tau)$ of finite index. Then there exists a \tb \ group topology $\tau'$ on $G$ \st :
	\begin{enumerate}
		\item \label{Tarski-like1} $H$ is a closed subgroup of $(G,\tau')$.
		\item \label{Tarski-like2} If $X$ is the character group of $(G,\tau)$, and $Y$ is the character group of $(G,\tau')$, then $Y$ is isomorphic to Hom$(G/H,\TT) \times X$.
		\item \label{Tarski-like3} The completion $K$ of $(G,\tau')$ contains a copy of $(G,\tau)$ of  index $[G:H]$.
		\item \label{Tarski-like4} There is a continuous epimorphism $k:K\to (G,\tau)$, extending the identity from $(G,\tau')$ onto $(G,\tau)$.
		\item \label{Tarski-like5} The \tg s $H$ as a subgroup of $K$, and $H$ as a subgroup of $(G,\tau)$ are the same.
	\end{enumerate}
\end{thm}

\begin{rem} \label{Remark3.17}  Notice that this result, which is essentially \cite[(4.15 ii)]{comf-hbook}, is reminiscent of Tarski's Paradox, obtaining from one compact Hausdorff group, $(G,\tau)$, a bigger group, but no much bigger, $K$, in which a copy of $G$ is dense in $K$, but at the same time $K$ is the finite union of cosets of $(G,\tau)$.
We offer a different proof than the one in \cite[(4.15 ii)]{comf-hbook}.
\end{rem}

\noindent {\em Proof of Theorem \ref{Tarski-like}.} Obviously Hom$(G/H,\TT)$ is finite. Let $\pi:G\to G/H$ be the natural map. Consider the finite subgroup of Hom$(G,\TT)$ given by
\[\sF:=\{f \circ \pi: f \in \mbox{Hom}(G/H,\TT)\}.\]
{Notice that $\psi \in \sF$ implies $\psi[H]=\{0\}$,} and since $H$ is a dense subgroup of $(G,\tau)$, it follows from Theorem \ref{B} that $\sF\cap X=\{0\}$. Set
\[Y:=\sF+X.\]
Finally, set $\tau':=\tau_Y$. Then $Y$ is the character of $(G,\tau')$ which is algebraically isomorphic to the character group of its compact Weil completion $K$.
Since $G/H$ is finite, it is immediate to see  that $H=\cap \{\ker f:f \in \sF\}$, hence $H$ is a closed subgroup of $(G,\tau')$, proving (\ref{Tarski-like1}). Since $\sF\cap X=\{0\}$ and $\sF$ is (clearly) isomorphic to {Hom}$(G/H,\TT)$, (\ref{Tarski-like2}) follows. Applying Pontryagin-van Kampen duality to (\ref{Tarski-like2}), it follows that $K$ is isomorphic to $G/H \times (G,\tau)$, yielding (\ref{Tarski-like3}), and (\ref{Tarski-like4}) in turn. (\ref{Tarski-like5}) follows {since $\psi \in \sF$ implies $\psi[H]=\{0\}$}. \epf

\bexm\label{3.11}  At the end of \cite{saxl-wil97}, the authors sketch a proof of the fact that, for any finite group $F$ and every non-principal ultrafilter $\sU$ of $\omega$, the (compact Hausdorff) group $G:=F^\omega$ contains a dense subgroup $G_\sU$ \st \  $G/G_\sU$ is isomorphic to $F$, hence  $G_\sU$ is of finite index in $G$. It follows from \cite{pospisil37a} that $G$ has $2^\cc$-many dense subgroups of finite index.
\eexm

\section{Closure of subgroups}

Our objective in this section is twofold:

\begin{enumerate}
\item [(a)] Let $G_S$ {be a precompact Abelian group with character group $S$. Characterize the} closure of a given subgroup of $G$ in terms of $S$.
\item [(b)] Let $\sF$ be a family of subgroups of $G$. Find the number of totally bounded group topologies $\tau$ on $G$ such that every element
of $\sF$ is closed in $(G,\tau)$.
\end{enumerate}

We start with a simple fact.

\begin{lem}\label{A1} A subgroup $H$ of an Abelian group $G$ is closed in $G_S$ \sii \ for all $a \in G\setminus H$ \ there is $\phi \in \AA(S,H)$ \st \ $\phi(a)\neq 0$.
\end{lem}

\pf \ $(\Leftarrow)$. Obvious. To see ($\Rightarrow$), notice that $G_S/H$ is Hausdorff \cite[(5.21)]{hewitt1963abstract} and precompact, hence \tb . By the Comfort-Ross Theorem, the character group of $G_S/H$ separates points. Therefore, there is a continuous $f:G_S/H\to \TT$ with $f(aH)\neq 0$. Then $\phi: g \mapsto f(gH)$ is as required. \epf.

Notice we are not assuming that $S$ is point-separating.
\medskip
\begin{cor}\label{3.1}  {Let $H\leq G$ and $S\leq\wG$. Then
		the following assertions are equivalent:
		\begin{enumerate}
			\item $g\in\overline{H}^{\,G_S}$.
			\item $\varphi(g)=0$ for all $\varphi\in\AA(S,H)$.
			\item $\varphi(g)\in\overline{\varphi [H]}^{\,\TT}$, whenever $\varphi\in S$.
		\end{enumerate}
	}
\end{cor}
\pf\,  (1) $\Rightarrow$ (2), (1) $\Rightarrow$ (3) and (3) $\Rightarrow$ (2) are obvious. On the other hand,
(2) $\Rightarrow$ (1) is a straightforward consequence of Lemma \ref{A1}.
\epf
\bkp

The next results are very similar to \cite[(23.24(c)) and (24.10)]{hewitt1963abstract} with very similar proofs.

\bthm\label{3.7} Let $S$ be a
subgroup of $\wG$, and $N:=\AA(G,S)$.
\begin{enumerate}
\item If $T$ is a subgroup of $S$, then $\AA(G,T)$ is a closed subgroup of $G_S$. In particular, $N=\overline{\{0\}}^{G_S}$.
\item { If $H\leq G$, then $\overline{H}^{G_S}=\AA(G,\AA(S,H))$. }
\item { If $H$ is a closed subgroup of $G_S$, then $H=H+N=\AA(G,\AA(S,H))$. }
\end{enumerate} In particular, if $H$ is a subgroup of $G$, then $H=\AA(G,\AA(\wG,H))$.
\ethm

\pf \ (1) If $x \in G\setminus \AA(G,T)$, then there is $\varphi \in T$ \st \ $\varphi(x)\neq 0$. Then apply Lemma \ref{A1}.

\noindent (2) ($\subseteq$) Let $g$ and $\phi$ be arbitrary elements of $\overline{H}^{G_S}$ and $\AA(S,H)$,
res\-pec\-ti\-ve\-ly.
Then $\phi(g)\in \overline{\phi[H]}=\{0\}$, which implies $g\in \AA(G,\AA(S,H))$.

($\supseteq$) Suppose that $g\notin \overline{H}^{G_S}$. Set $K:={\overline{G_S}}$, the Weil completion of $G_S$,
which is a compact group, and $L:=\overline{H}^{K}$. We have that $g\in K\setminus L$ and $K\: \^ =S$.
Therefore, by \cite[(23.26)]{hewitt1963abstract} there is $\chi\in K\: \^ \simeq S$ \st\ $\chi[L]=\{0\}$ and $\chi(g)\not= 0$. Plainly, $\chi\in \AA(S,H)$
and $g\notin  \AA(G,\AA(S,H))$.

\noindent (3) {follows from (2).}

The last statement follows by taking $S=\wG$ since in this case all subgroups of $G$ are closed in $G_S$ (Lemma \ref{CoSa}). \epf

\brem\label{4.11}
One may be tempted to believe that if $S$ is a point-separating subgroup of $\wG$, then $S=\AA(\wG,\AA(G,S))$. This is false in general.
Of course $S\subseteq \AA(\wG,\AA(G,S))$. But if $S$ is any non-torsion subgroup of $\TT \simeq \wZ$, then $\AA(\ZZ,S)=\{0\} \implies \AA(\wZ,\AA(\ZZ,S))=\wZ$; see Corollary \ref{4.3} below.
\erem

\begin{cor}\label{3.8.1}  {Let $H\leq N$ be subgroups of $G$, and let $S$ be a 
 subgroup of $\wG$. }
		The following conditions are equivalent:
		\begin{enumerate}
			\item $N\subseteq \overline H^{G_S}$.
			\item $H$ is $\tau_S$-dense in $N$.
			\item $\AA(G,\AA(S,N))=\AA(G,\AA(S,H))$.
		\end{enumerate}

\noindent {In addition, $N=\overline H^{G_S}$ if and only if $N=\AA(G,\AA(S,N))=\AA(G,\AA(S,H))$.}\epf
\end{cor}
\mkp

Comfort and Saks were the first who gave examples of SC groups.
\blem\label{CoSa}{\rm(\cite[Lemma 2.1]{comfsaks73})}
Let $G$ be an Abelian group with the Bohr topology $\tau_b$. Then every subgroup of
$G$ is closed in $(G,\tau_b)$.
\elem

By using this result the following is shown in \cite[Lemma 2.5]{CDT}: Let $G$ be an Abelian group of bounded order, and let $\tau$ be a \tb\ group
topology on $G$. Then every subgroup of $G$ is closed in $(G,\tau)$ if and only if $\tau=\tau_b$.\noindent

Hence for such groups $G$ the following is false: Let $\{H_i: i\in I\}$ be a family of subgroups of $G$. Then there is a \tb\ group topology
$\tau\neq\tau_b$ on $G$ such that $H_i$ is closed in $(G,\tau)$ for all $i\in I$.\noindent

In Theorem \ref{Thm_A} we give the best possible result for $|I|<2^{|G|}$. The corresponding proof needs some preparation.

\blem\label{CO}{\rm(\cite[Theorem 4.3.]{BerhanuComfortReid1985})}
Let $(G,\tau)$ be an infinite \tb \ group. Then
\noindent
\begin{enumerate}
\item [(a)] $\chi(G,\tau)=w(G,\tau)$.
\item [(b)] If in addition $G$ is Abelian, then $\chi(G,\tau)= (G,\tau) \^ $.
\end{enumerate}
\elem

\blem\label{D1}{\rm(\cite[Lemma (2.9)]{Remus1991})}
Let $G$ be a discrete group which is maximally almost periodic. Then $G$
admits a totally bounded group topology $\tau$ with $w(G,\tau)\leq|G|$.
\elem

\blem\label{D2} {\rm(\cite[Lemma (2.16)]{Remus1983})}
Let $G$ be an infinite Abelian group. If $r(G)$ is finite, then $G$ is countable.
Otherwise, $r(G)=|G|$ holds.
\elem

\blem\label{SUP}
Let $G$ be an infinite group, and let $\{\tau_i:i\in I\}$ be a family of group topologies on $G$ such there is $i_0\in I$ with $\chi(G,\tau_{i_0})\geq\omega$. Set $\tau:=\bigvee_{i\in I}\tau_i$ in the lattice of all group topologies on $G$.
Then $\chi(G,\tau)\leq\max\{|I|,\sup\{\chi(G,\tau_i):i\in I\}\}$.
\elem
\pf\
Let $e$ be the identity of $G$. For all $i\in I$, let $\sB_i$ be a \nhd-basis of $e$ with respect to $(G,\tau_i)$ such that $|\sB_i|=\chi(G,\tau_i)$. For any finite subset $M=\{i_1,\ldots,i_n\}$
of $I$ let $\sB_M:=\{\bigcap^{i_n}_{j=i_1}U_j:U_j\in \sB_j\}$. Let $\sF$  be the set of all finite subsets of $I$.
Then $\sB:=\bigcup_{M\in\sF}\sB_M$ is a \nhd-basis of $e$ with respect to $\tau$. Let $\alpha:=\sup\{\chi(G,\tau_i):i\in I\}$.
Then $\alpha\geq\omega$, $|\sB_M|\leq\alpha$ for all $M\in\sF$, and $|\sB_{\{i_0\}}|=|\sB_{i_0}|\geq\omega$. For $|I|\geq\omega$ we have $|\sF|=|I|$.
Since $\sB:=\bigcup_{M\in\sF}\sB_M$, finally \cite[Proposition 4.4 and 4.5, Chapter III]{Levy1979} implies
$|\sB|\leq\max\{|\sF|,\alpha\}=\max\{|I|,\alpha\}$.
\epf
\bkp

{We are ready for the proof of our first main result}.
\bkp

\noindent \textbf{Proof of Theorem \ref{Thm_A}}.
Let $i\in I$. By Lemma \ref{D1}, $G$ and $G/H_i$ admit  totally bounded group topologies $\tau$ and $\nu_i$, respectively, such that $w(G,\tau)\leq|G|$ and
$w(G/H,\nu_i)\leq|G/H_I|\leq|G|$. Let $\pi_i:G\rightarrow G/H_i$ be the canonical epimorphism and $\tilde{\nu_i}$ be the initial topology
on $G$ \wrt \ $\pi_i$ and $\nu_i$. Surely $\tilde{\nu_i}$ is precompact.  Then the quotient topology of $\tilde{\nu_i}$ on $G/H_i$ is $\nu_i$. Hence $H_i$ is closed in $(G,\tilde{\nu_i})$.
Now $\chi(G,\tilde{\nu_i})=\chi(G/H,\nu_i)\leq|G|$. Define $\mu_0:=\tau\vee(\bigvee_{i\in I}\tilde{\nu_i})$. Then $\mu_0$ is totally bounded, and $H_i$ is closed
in $(G,\mu_0)$. Lemma \ref{SUP} implies $\chi(G,\mu_0)\leq\max\{|I|,|G|\}$. By $|I|<2^{|G|}$ we get $\chi(G,\mu_0)<2^{|G|}$.
Finally, Lemma \ref{CO} implies $w(G,\mu_0)=|(G,\mu_0)\sphat\,|<2^{|G|}$.

$|\wG|=2^{|G|}$ holds by a result of Kakutani \cite{kaku43}. For $M:=(G,\mu_0)\sphat$
we have $|\wG/M|=2^{|G|}$ since $|M|<2^{|G|}$. Lemma~\ref{D2} implies $r(\wG/M)=2^{|G|}$. Thus $\wG/M$ contains $2^{2^{|G|}}$-many subgroups.
Hence $\wG$ has this number of subgroups containing $M$.  Thus by the Comfort-Ross Theorem there are $2^{2^{|G|}}$-many totally bounded
group topologies $\mu$ on $G$ being finer then $\mu_0$. Surely $H_i$ is closed in $(G,\mu)$.
\epf
\mkp

The argumentation in the first part of the proof of Theorem~\ref{Thm_A} (for $|I|=1$) shows:

\begin{thm}\label{D4}
{Let $G$ be an infinite group, and let $H$ be a normal subgroup of $G$ such that $G/H$ is maximally almost
periodic in the discrete topology. If the finest precompact group topology on $G$ is Hausdorff,
then there is a totally bounded group topology $\mu_0$ on $G$ with $w(G,\mu_0)\leq|G|$ such that $H$ is closed in $(G,\mu_0)$}.\epf
\end{thm}

By applying Theorem \ref{D4} and \cite[Corollary 2.6(a)]{ComfortRemus1996}, we get

\begin{thm}\label{D5}
{Let $G$ be an infinite group, and let $H$ be a normal subgroup of $G$ such that $G/H$ is maximally almost
periodic in the discrete topology.
If the finest precompact group topology $\tau_f$ on $G$ is Hausdorff with $\alpha :=w(G,\tau_f)>|G|$, then $G$ admits at least
$\alpha$-many totally group topologies $\mu$ such that $H$ is closed in $(G,\mu)$.}\epf
\end{thm}

\brem\label{D6}
Remus showed in \cite[Satz (2.2)]{RemusNovember1995} that if $G$ is an Abelian group and $H$ is a subgroup of $G$, then for every $\tau\in\sP\sK(H)$ there is a finest precompact group topology
$\tau_f^H$ on $G$ which induces $\tau$. If $\tau$ is totally bounded, then $\tau_f^H$ has the same property. This implies that all subgroups of $G$ which contain $H$ are closed in $(G,\tau_f^H)$, see \cite[Folgerung (2.6)(a)]{RemusNovember1995}.

In summary: If $G$ is an Abelian group and $H$ is an infinite subgroup of $G$, then there are
at least $2^{2^{|H|}}$-many totally bounded group topologies $\tau$ on $G$ such that every subgroup $N$ of $G$ with $H\subseteq N$ is closed in $(G,\tau)$.
\erem
\bkp

{Let $G$ be an infinite countable Abelian group with only countably many subgroups. These groups are classified in \cite{Rychkov1991}.
Groups like $\bigoplus_m\ZZ\oplus\bigoplus_nZ(p^\infty)$ with $m,n<\omega$ are examples. By Theorem \ref{Thm_A} these groups have exactly $2^{\cc}$ many SC-group topologies.}

\section{Subgroups of $\wG$ producing the same closed subgroups in $G$}

{In this section, we study the poset of all subgroups $S\leq \wG$ having a pre\-vious\-ly fixed subgroup $H\leq G$ closed in the weak topology $\tau_S$.
We prove that there always exists a greatest element (Theorem \ref{Thm_B}) but not a smallest one in general.}\bkp

{We first characterize the subgroups of $\wG$ producing the same closed subgroups of $G$.}\bkp

\bthm\label{5.3}
{Let $G$ be an Abelian group and let $S_1, S_2$ be subgroups of $\wG$. Then $S_1$ and $S_2$ produce the same closed subgroups in $G$
	if and only if $$\overline{L\cap S_1}^{\,\wG}=\overline{L\cap S_2}^{\,\wG}$$
	for all closed subgroups $L$ of $\wG$, when equipped with $\gs(\wG,G)$.}
\ethm
\pf \, ($\Leftarrow$) Suppose that $$\overline{L\cap S_1}^{\,\wG}=\overline{L\cap S_2}^{\,\wG}$$
for all closed subgroups $L$ of $\wG$ and let $H$ be a $\tau_{S_1}$-closed subgroup of $G$.
Then $S_1\cap \AA(\wG,H) = \AA(S_1,H)$ is dense in $\AA(\wG,H)$ by Lemma \ref{AA}. As a consequence,
$S_2\cap \AA(\wG,H) = \AA(S_2,H)$ is dense in $\AA(\wG,H)$ as well, which implies by Lemma \ref{AA} that $H$ is $\tau_{S_2}$-closed.

\noindent ($\Rightarrow$)  Reasoning by contradiction, suppose there is a closed subgroup $L$ of $\wG$ such that
$\overline{S_1\cap L}^{\, \wG}\nsubseteq \overline{S_2\cap L}^{\, \wG}$. Observe that, by Lemma \ref{closure},
$$\overline{S_i\cap L}^{\, \wG}=\AA(\wG,\AA(G, \overline{S_i\cap L}^{\, \wG})),\ 1\leq i\leq 2.$$

\noindent Set $H:=\AA(G, \overline{S_1\cap L}^{\, \wG})$.
Then $H$ is $\tau_{S_1}$-closed by Theorem \ref{3.7}. However, Lemma \ref{AA} implies
$$\overline{\AA(S_2,H)}^{\,\wG}=\overline{S_2\cap \AA(\wG,H)}^{\,\wG}=
\overline{S_2\cap \AA(\wG,\AA(G, \overline{S_1\cap L}^{\, \wG}))}^{\,\wG}=
\overline{S_2\cap \overline{S_1\cap L}^{\, \wG}}^{\,\wG}\subseteq$$
$$\overline{S_2\cap L}^{\,\wG}\nsupseteq \overline{S_1\cap L}^{\,\wG}=\AA(\wG,\AA(G, \overline{S_1\cap L}^{\, \wG}))=
\AA(\wG,H).$$

As a consequence, $\AA(S_2,H)$ is not dense in $\AA(\wG,H)$, which means by Lemma \ref{AA} that $H$ is not $\tau_{S_2}$-closed.
This completes the proof.
\epf
\bkp

{The latter result implies the following consequences.}
In the Corollary below, note that $G/\AA(G,R)$ is Hausdorff.

\bcor\label{5.2.3}
{Let $S\leq \wG$ and let $R:= \overline{S}^{\wG}$. The following assertions are equivalent:
	\begin{enumerate}
		\item $\overline{L\cap S}^{\,\wG}=L\cap R$ for every closed subgroup $L$ of $\wG$.
		\item $S$ and $R$ produce the same closed subgroups in $G$.
		\item $S$ and $R$ produce the same closed subgroups in $G/\AA(G,R)$.
	\end{enumerate}
}
\ecor
\pf\, \noindent (1)$\Leftrightarrow $ (2) is a straigtforward consequence of Theorem \ref{5.3}.

\noindent (2)$\Leftrightarrow $ (3) By Theorem \ref{3.7} the closure of $\{0\}$ in $G_S$ resp. $G_R$ is $\AA(G,R)$.
Now apply \cite[(5.34)]{hewitt1963abstract} .
\epf

\bexm\label{5.3Example}  Fix different prime numbers $p$ and $q$ and consider $S:=Z(p^\infty)$. Of course, $R:=\overline{S}^\wZ=\TT$. We will see below in Corollary \ref{4.2} (or see \cite[(3.5.4) and (3.5.5)]{Dik-Prod-Stoy1990}), that the closed subgroups of $\ZZ_S$ have the form $p^k\ZZ$. On the other hand, every subgroup of $\ZZ$ is closed in $\ZZ_R$ (Lemma \ref{CoSa}. If $L:=\< 1/q\> $, then $L$ is a closed subgroup of $\wZ$, but $\overline{L\cap S}^{\,\wZ}=\{0\}\neq L=L\cap R$.
\eexm

\bcor\label{5.2.1}
{Let $R\leq \wG$ be a closed subgroup of $\wG$. Then $H\leq G$ is $\tau_R$-closed if and only if
	$H=H+\AA(G,R)$. In other words, $\overline{H}^{\,G_R}=H+\AA(G,R)$.}
\ecor
\pf\, If $R$ separates points in $G$, then $R=\wG$ and the result follows from \cite[Lemma 2.1]{comfsaks73}. So we may assume that $R\not= \wG$.
By Theorem {\ref{3.7},} $\overline{\{0\}}^{G_R}=\AA(G,R)$.
Take the canonical continuous epimorphism $\pi:(G,\tau_R)\rightarrow (G/\AA(G,R),\tau_R^q)$, where $\tau_R^q$ is the Hausdorff quotient topology.
Now $\AA(G,R)\subseteq \overline{H}^{G_R}$ implies
$H\subseteq H+\AA(G,R)\subseteq \overline{H}^{\,G_R}$. By \cite[(23.25)]{hewitt1963abstract} the character group of $G/A(G,R)$ endowed with the
discrete topology is topologically isomorphic to $\AA(\wG,\AA(G,R))$. Lemma \ref{closure} implies $R=\AA(\wG,\AA(G,R))$.
The topological isomorphism is definded by $\rho(\psi)=\psi\circ\pi$. Thus $\tau_R^q$ is the Bohr topology on $G/\AA(G,R)$.
As a consequence, {every} subgroup
of the latter group is $\tau_R$-closed by Lemma \ref{CoSa}. In particular, the subgroup $\pi[H]$ is $\tau_R$-closed. Being
the quotient map obviously $\tau_R$-continuous, it follows that
$$H+\AA(G,R)=\pi^{-1}[\pi[H]]$$ is $\tau_R$-closed in $G$.
\epf
\bkp

{We can now use Theorem \ref{5.3} as a main tool to prove that every subgroup $S\leq\wG$ is contained in a greatest subgroup
$MS$ defining the same set of	$\tau_S$-closed subgroups.}
\bkp

\noindent \textbf{Proof of Theorem \ref{Thm_B}}.
Set $$\mathcal{S}:=\{T\leq\wG : S\leq T\ \hbox{and}\ C_T=C_S\}.$$
We claim that the pair $(\mathcal{S},\subseteq )$ is inductive if it is ordered by inclusion.

Indeed, let $\{S_i: i\in I\}$ be a chain in $(\mathcal{S},\subseteq )$. Take $S_0:=\cup\{S_i : i\in I\}$.
If $L\leq\wG$ is closed, we have $$\overline{S\cap L}^{\, \wG}\subseteq \overline{S_0\cap L}^{\, \wG}.$$
On the other hand
$$\overline{S_0\cap L}^{\, \wG}=\overline{(\cup\{S_i : i\in I\})\cap L}^{\, \wG}=
\overline{\cup\{S_i\cap L : i\in I\}}^{\, \wG}\subseteq$$
$$\overline{\cup\{\overline{S_i\cap L} : i\in I\}}^{\, \wG}=\overline{\cup\{\overline{S\cap L} : i\in I\}}^{\, \wG}=
\overline{S\cap L}{\, \wG}.$$
{The second to last equality above follows from Theorem \ref{5.3}.} This implies that $S_0\in\mathcal{S}$, which completes the proof of the claim.
Therefore, we have verified the existence of maximal elements in $(\mathcal{S},\subseteq )$ by Zorn's Lemma.
In order to demonstrate that there is a greatest element, it will suffice to prove that
if $S_1$ and $S_2$ are in $\mathcal S$, so is $S_1+S_2$.

It is clear that every subgroup $\tau_S$-closed is also $\tau_{(S_1+S_2)}$-closed. Therefore
$C_S\subseteq C_{(S_1+S_2)}$.

Conversely, suppose that $H\leq G$ is $\tau_{(S_1+S_2)}$-closed. By Lemma \ref{A1}
we have that if $g\notin H$, there is $\chi_1\in S_1$ and $\chi_2\in S_2$ such that $(\chi_1+\chi_2)[H]=\{0\}$ and
$(\chi_1+\chi_2)(g)\not=0$. Therefore $$\chi_{1|H}=-\chi_{2|H}.$$
Set $$K=\overline{H}^{G_S}$$
and observe that, since $S_i\in\mathcal{S}$, $1\leq i\leq 2$,
it follows
$$K=\overline{H}^{G_S}=\overline{H}^{G_{S_i}},\ 1\leq i\leq 2.$$
Therefore both characters $\chi_1$ and $\chi_2$ can be extended to $K$ continuously.
Denoting these extensions by $\chi_1$ and $\chi_2$ again for simplicity's sake, we have that
$(\chi_1+\chi_2)[H]=\{0\}$, which implies $(\chi_1+\chi_2)[K]=\{0\}$.
This entails that $g\notin K$. Hence $K=H$, which completes the proof.
\epf
\bkp

\brem\label{Remark 5.7}
{Given an Abelian group $G$, it is not true in general that for each subgroup $S\leq \wG$,
there exists a minimum subgroup $mS\leq \wG$ such that $$C_{mS}=C_S$$ (see Theorem \ref{4.7} below).}
\erem

\section{Totally bounded groups topologies in which every subgroup is closed}

In this section we approach the following

\bqtn\label{Q1}
Let $G$ be an infinite Abelian group.
\begin{enumerate}
\item [(a)] Characterize the subgroups $S$ of $\wG$ \st \
all subgroups of $G$ are $\tau_{S}$-closed.
\item [(b)] Find the number of SC-group topologies on $G$.
\end{enumerate}
\eqtn
\mkp

Consider the torsion subgroup $t\wG$ of $\wG$, take $\phi \in t\wG$ and $H:=$~ker~$\phi$. We have then that $H$ is of finite index in $G$,
hence if it were closed in $G_{S}$, then it would be also open and this would imply $\phi \in S$. It follows

\blem\label{Lemma 1.}
{If all the subgroups of $G$ are $\tau_{S}$-closed, then $t\wG \leq S$.}\epf
\elem
\mkp

As a consequence, we have

\bcor\label{Corollary 1.}
{ If $G$ is of bounded order and all the subgroups of $G$ are $\tau_{S}$-closed, then $S=\wG$ and $G_S=(G,\tau_b(G))$.}\epf
\ecor
\mkp

After this result, one might conjecture that {all the subgroups of $G$ are $\tau_{S}$-closed \sii \  $t\wG \leq S$.}
However this is wrong. Indeed, consider the group $G=Z(p^\infty), p$ prime. Since all proper subgroups of $G$ are finite,
we have that all subgroups are also closed on every totally bounded group topology of $G$. Nevertheless, we have that
$\wG=\Delta_p$, which is a torsion-free group. This shows that $t\wG$ can be trivial even though every totally bounded topology on $G$
has all its subgroups closed.

\bdfn\label{Def_5.3} (\cite{Dik-Prod-Stoy1990}, p.~133)
A subgroup $H$ of a topological group $(G,\tau)$ is \emph{totally dense} in $G$ if $H\cap K$ is dense in $K$ for every closed normal subgroup $K$ of $G$.
\edfn

Now Question \ref{Q1}(a) can be solved.

\bprp\label{5.3.1}
Let $G$ be an Abelian group and $S$ be a subgroup of $\widehat G$.
Then every subgroup of $(G,\tau_S)$ is closed if and only if $S$ is totally dense in $\widehat G$.
\eprp
\pf\ {By Lemma \ref{CoSa} every subgroup of $G$ is closed in $(G,\tau_{\hat{G}})$.
Then apply Theorem \ref{5.3} if we take $S_1=S$ and $S_2=\wG$}.\epf
\mkp

Now \cite[Corollary 2]{Khan1984} and Propositon \ref{5.3.1} imply

\bprp\label{khan1}
Let $G$ be an Abelian group equipped with the finest precompact group topology, that is, the Bohr topology $\tau_b(G)$.
Then the following assertions are equivalent:
\begin{enumerate}
	\item $\tau_b(G)$ is the only precompact group topology $\tau$ on $G$ such that all subgroups of $(G,\tau)$ are closed;
	\item $G$ is of bounded order.
\end{enumerate}
\eprp

\mkp

We notice that  Corollary \ref{Corollary 1.} yields the implication $(2)\Rightarrow (1)$ in the above proposition that,
incidentally, is proven in \cite[Lemma 2.5]{CDT}.
On the other hand, Proposition \ref{5.3.1} and \cite[Exercise 5.5.6]{Dik-Prod-Stoy1990} imply

\bprp\label{khan2}
Let $G$ be an infinite Abelian group. Then the following assertions are equivalent:
\begin{enumerate}
	\item For every totally bounded group $(G,\tau)$ each subgroup of $(G,\tau)$ is closed;
	\item There is a prime number $p$ such that $G$ is isomorphic to $Z(p^\infty)$.
\end{enumerate}
\eprp

\mkp

From Proposition \ref{5.3.1}, the solution to Question \ref{Q1}(b) for a given group $G$ reduces to the search of totally dense subgroups of $\wG$.
In this direction, Comfort and Dikranjan \cite{CD1999} have proven that the smallest (under inclusion) totally dense subgroup of a compact Abelian group
$K$ is $tK$, {\em when it is itself totally dense.} Furthermore, according to their Theorem 4.1, this happens \sii \ $K$ has no copies of $\Delta_p$.
Therefore, we have the following:

\bprp\label{Proposition 1.}  {Let $G$ be an Abelian group \st \  $\wG$ has copies of $\Delta_p$ for no prime $p$. Let $S$ be a subgroup of $\wG$. All the subgroups of $G$ are $\tau_{S}$-closed \sii \ $t\wG \leq S \leq \wG$.}
\eprp
\mkp

 Every \tb \ group topology on $G=Z(p^\infty)$ has all its subgroups closed. In this case of course there exist $2^\cc$-many such topologies.  Hence we see that the number of \tb \ group topologies making all the subgroups of $G$ closed depends on $G$ and goes from 1 (groups of bounded order) all the way to $2^{2^{|G|}}$ ($G=Z(p^\infty)$).

In relation with the above questions, we first list some observations.

\brem\label{newrem 1}
\indent
\begin{enumerate}
\item[(a)] Let $G$ be an Abelian group which is not of bounded order. Then $G$ admits at least $\cc$-many \tb \ group topologies $\tau$ \st \ every subgroup of $G$ is closed in $(G,\tau).$ For, $\wG=B(\wG)$, the set of compact elements of $\wG$. If $\wG$ were not {\em admissible}, \ie \ if $B(\wG)=t\wG$, then $\wG$ would be a torsion group, implying that  $G$ is of bounded order \cite[(25.9)]{hewitt1963abstract}, a contradiction. By \cite[Corollary 3]{Khan1984} $\wG$ contains a totally dense subgroup $H$ such that $\wG/H \cong \QQ$. Since $\QQ$ contains $\cc$-many subgroups, we see that $\wG$ contains $\cc$-many subgroups, each containing the totally dense subgroup $H$. Now apply Proposition \ref{5.3.1}.

\item[(b)] In \cite[p.~170]{soundaiii} \textit{SC groups} (see Section 1) and \textit{DSC groups} are introduced: A compact Abelian group having a dense SC group is called a DSC group, in other words: it is the completion of a SC group.

Let $G$ be any group. Then the finite-index topology $\tau_l$ on $G$ is the finest linear precompact topology on $G$. If $G$ is Abelian, the character group of $(G,\tau_l)$ is $t\wG$. Proposition \ref{5.3.1} implies: $(G,\tau_l)$ is a SC group \sii \  $t\wG$ is totally dense in $\wG$. Note that $\tau_l$ is the only possible linear precompact group topology $\mu$ such that $(G,\mu)$ is a SC group.

Now the proof of \cite[Theorem 1.15]{soundaiii}  implies the following:
$(G,\tau_l)$ is a SC group \sii \ $G$ has a free Abelian subgroup $F$ of finite rank such that $G/F$ is a torsion group and for each prime $p$
the $p$-component of $G/F$ is of bounded order. For the proof, take $H=(G,\tau_l)$, and $G$ the completion of this $H$ in the proof of \cite[Theorem 1.15]{soundaiii}.
Note that in the proof of \cite[Theorem 1.15]{soundaiii} it is allowed that the rank of $F$ is zero
[Look at the end of the proof of \cite[Proposition 1.14]{soundaiii}.].
\end{enumerate}
\erem
\mkp

First we will show that Abelian groups $G$ with $r_0(G)>0$ have many SC-group topologies.

\blem\label{newlem 1}
Let $m>0$ be a cardinal. Then $\TT^m$ contains $(2^\cc\cdot~2^{2^m})$-many totally dense subgroups.
\elem
\pf\
By \cite[Theorem 3]{DikranyanStoyanov1986} or \cite{Grant1979}, the group $(t\TT)^m$ is totally dense in $\TT^m$.
Let $H:=\TT^m/(t\TT)^m$. Then $H$ is algebraically isomorphic to $\RR^m$ and, therefore, has torsion-free rank $|H|$.
Hence there are $2^{|H|}$-many subgroups $S$ of $\TT^m$ containing the totally dense group $(t\TT)^m$.
Thus each such $S$ is a totally dense subgroup of $\TT^m$. Finally, it suffices to notice that $|H|=\cc^m=\cc\cdot 2^m$.
\epf

\bthm\label{newthm 1}
Let $G$ be an Abelian group with $r_0(G)>0$. Then there are $(2^\cc\cdot 2^{2^{r_0(G)}})$-many SC-group topologies on $G$.
\ethm
\pf\ Let $K:=\hat{G}$. By \cite[Lemma 5.4]{ComfortRemus1994a}, there is a continuous homomorphism from $K$ onto $\TT^{\,r_0(G)}$.
Now it suffices to apply \cite[Lemma 4.1 (c)]{comfsound82}, Lemma \ref{newlem 1} and Proposition \ref{5.3.1} in order to complete the proof.
\epf

\bcor\label{newcor 1}
Let $G$ be an infinite Abelian group.
If $r_0(G)=|G|$, then there are exactly $2^{2^{|G|}}$ -many SC-group topologies on $G$.\epf
\ecor
\mkp

Next we consider some Abelian torsion groups that are not of bounded order.

\bthm\label{nnnewthm 1}
Let $G$ be an Abelian torsion group such that $I=\{p\in\PP:r_p(G)\neq 0\}$ is infinite.
Then $G$ admits at least $2^\cc$-many SC-group topologies.
\ethm

\pf\
Let $H:=\prod_{p\in I}Z(p)$. By \cite[Lemma 5.5]{comfsound82}, $tH=\bigoplus_{p\in I}Z(p)$ is totally dense in $H$. Since $|tH|=\omega<|H|=\cc$,
the group $H/tH$ has $2^\cc$ many subgroups by Lemma 4.8. Hence $H$ possesses the same number of totally dense subgroups.
Let $K:=\widehat{G}$.
Then by \cite [Lemma 5.4]{ComfortRemus1994a} there is a continuous epimorphism from $K$ onto $\prod_{p\in\PP}(Z(p))^{r_p(G)}$. In particular, there is a continuous
epimorphism from $K$ onto $H$. Now \cite[Lemma 4.1(c)]{comfsound82} implies that $K$ has $2^\cc$-many totally dense subgroups. Finally apply Proposition \ref{5.3.1} to complete the proof.
\epf
\bkp

We observe that Proposition \ref{khan1} implies that the above result fails if $I$ is finite.

\bprp\label{newthm 2}
Let $G$ be an infinite Abelian torsion group of cardinality $\alpha$. Assume {\rm cf}$(\alpha)=\omega$ and $\log(2^\alpha)=\alpha$.
Let $I:=\{p\in\PP:r_p(G)\geq\omega\}$ be infinite. If $r_p(G)< |G|$ for all $p\in I$ and all $r_p(G)$ with $p\in I$ are distinct, then $G$ admits exactly $2^{2^{|G|}}$-many
SC-group topologies.
\eprp
\pf\
Let $\alpha_p:=r_p(G)$ for all $p\in\PP$. Consider $K:=\prod_{p\in I}Z(p)^{\alpha_p}$. Then $tK=\bigoplus_{p\in I}(\bigoplus_{2^{\alpha_p}}Z(p))$.
Hence $\beta:=|tK|=\sum_{p\in I}2^{\alpha_p}$. Since $G$ has infinite rank, by Lemma \ref{D2} we get $\alpha=\sum_{p\in I}\alpha_p$.
Then \cite[Lemma 5.4]{ComfortRobertson1985} implies $\beta <2^\alpha$. $K$ is the character group of $H:=\bigoplus_{p\in I}(\bigoplus_{\alpha_p}Z(p))$.
Thus $|K|=2^\alpha$. Hence $|tK|<|K|$. By \cite[Lemma 5.5]{comfsound82}, $tK$ is totally dense in $K$.
Since $|K/tK|=|K|$, it follows
that $K$ has $2^{|K|}=2^{2^{|G|}}$ many totally dense subgroups. By \cite[Lemma 5.4]{ComfortRemus1994a}, there is a continuous epimorphism from  $G$ onto $K$.
Apply \cite[Lemma 4.1(c)]{comfsound82} and Proposition \ref{5.3.1} to complete the proof.
\epf

\brem\label{nnewrem 1}
(a) We give an example for an application of Theorem \ref{newthm 2}:

Let $m$ be an infinite cardinal. Define $\alpha_0:=m$ and $\alpha_{n+1}:=2^{\alpha_n}$ for all $n\in\omega$. Let $(p_n)_{n\in\omega}$
be a sequence of prime numbers which are pairwise different.
Set $G:=\bigoplus_{n\in\omega}(\bigoplus_{\alpha_n}Z(p_n))$. Then $|G|=\alpha:=\sum_{n\in\omega}\alpha_n>m$. Clearly cf$(\alpha)=\omega$,
and \cite[Lemma 5.4]{ComfortRobertson1985} implies $\alpha=\log(2^\alpha)$. Then Proposition \ref{newthm 2} implies that $G$ admits exactly $2^{2^{|G|}}$-many
SC-group topologies.

(b) Let $K$ be the Abelian compact totally disconnected group of weight $\alpha$ defined in the proof of the above Theorem.
Then cf$(\alpha)=\omega$ and $\alpha=\log(2^\alpha)$ imply the crucial fact $|tK|<|K|$.
Now let $L$ be any Abelian compact totally disconnected group of infinite weight $\alpha$ with $|tL|<|L|$. Then \cite[Theorem 5.8]{ComfortRobertson1985}
implies cf$(\alpha)=\omega$ and $\alpha=\log(2^\alpha)$.
\erem
\mkp

Now it is natural to pose the following

\bqtn\label{Q2}
Let $G$ be an infinite Abelian group which is not of bounded order. Does $G$ admit exactly $2^{2^{|G|}}$-many SC-group topologies in the
following cases:
\begin{enumerate}
\item [(a)] $G$ is a torsion {group,
\item [(b)] $0<r_0(G)<|G|$.}
\end{enumerate}
\eqtn
\mkp

We will show that in both cases the answer is ``no''. For that we need several lemmas.

\blem\label{Prufer}
For an infinite cardinal $m$ and $p\in\PP$ let $G:=\bigoplus_m Z(p^\infty)$. Then $G$ admits exactly 
$2^{2^m}$-many SC-group topologies.
\elem

\pf\
We have $\wG=\Delta_p^m$. By \cite[Lemma 5]{Khan1984} $\Delta_p$ contains a dense subgroup $H$ with $\Delta_p/H\cong\QQ$.
Then $H^m$ is dense in $\wG$, and $V:=\wG/H^m\cong\QQ^m\cong\bigoplus_{2^m}\QQ$. Let $$f\colon \Delta_p^m\to \bigoplus_{2^m}\QQ$$
be the algebraic epimorphism $$\Delta_p^m\to V\to \QQ^m\to \bigoplus_{2^m}\QQ.$$
If $S\subseteq 2^m$ is nonempty, let $\pi_S\colon \bigoplus_{2^m}\QQ\to \bigoplus_{S}\QQ$ be the natural projection and
let $\Sigma_S\colon \bigoplus_{S}\QQ\to \QQ$ be the homomorphism that adds up the coordinates of elements in $\bigoplus_{S}\QQ$:
$$\Sigma_S(q_j)=\sum\limits_{j\in S} q_j<\infty.$$
We then have
\begin{equation*}
\begin{picture}(148,90)
\put(-100,60){$\Delta_p^m$}
\put(-35,60){$\bigoplus_{2^m} \QQ$}
\put(-80,62){\vector(1,0){40}}
\put(-65,70){$f$}
\put(50,60){$\bigoplus_{S} \QQ$}
\put(05,62){\vector(1,0){40}}
\put(20,70){$\pi_S$}
\put(135,60){$\QQ.$}
\put(90,62){\vector(1,0){40}}
\put(102,70){$\Sigma_S$}
\end{picture}
\end{equation*}\vspace{-2cm}

Let $H_S := \ker (\Sigma_S\circ \pi_S\circ f)$. Note that $H^m\subseteq H_S$, hence $H_S$ is dense in $\Delta_p^m$.  \cite[Cor. 3]{Khan1984} implies that $H_S$ is totally dense in $\Delta_p^m$, since $\Delta_p^m/H_S\cong \QQ$.
We now show that if $S_1$ and $S_2$ are different subsets of $2^m$, then $H_{S_1}\not= H_{S_2}$.
Suppose that $\sigma\in S_1\setminus S_2$. Consider $q\in\QQ$, $q\not=0$,
and define $g=(g(t))\in \bigoplus_{2^m} \QQ$ by
\begin{equation*}
g(t):=\left\{
\begin{array}{ll}
q &\quad \hbox{if}\ t= \sigma, \\
0 &\quad \hbox{otherwise}.
\end{array}
\right.
\end{equation*}

Then we have that

\begin{equation*}
(\pi_{S_1}(g))(t)=\left\{
\begin{array}{ll}
q &\quad \hbox{if}\ t= \sigma, \\
0 &\quad \hbox{otherwise},
\end{array}
\right.
\end{equation*}

\noindent while  {$$(\pi_{S_2}(g))(t)=(0).$$}

Pick $x\in f^{-1}[\{g\}]$. Then $$(\Sigma_{S_1}\circ \pi_{S_1}\circ f)(x)=(\Sigma_{S_1}\circ\pi_{S_1})(g)=q\not= 0\Rightarrow x\notin H_{S_1},$$ while
{$$(\Sigma_{S_2}\circ \pi_{S_2}\circ f)(x)=(\Sigma_{S_2}\circ\pi_{S_2})(g) = \Sigma_{S_2}((0))= 0\Rightarrow x\in H_{S_2}.$$}\\
\noindent Hence $H_{S_1}\not=H_{S_2}$.

\noindent In summary, we have  constructed $2^{2^m}$-many totally dense subgroups $H_S$ of $\Delta_p^m$.
Now, applying Proposition \ref{5.3.1}, we obtain that $G$ admits exactly $2^{2^m}$-many SC-group topologies.
\epf
\mkp




\blem\label{nnewthm 2}
Let $G$ be an infinite Abelian group and $H$ an infinite subgroup of it.
If $H$ has $m$-many SC group topologies, then $G$ admits at least $m$-many such group topologies.
\elem

\pf\
By duality theory there is a continuous epimorphism from $\widehat{G}$ onto $\widehat{H}$.
Then Proposition 6.5 implies that $\widehat{H}$ has $m$-many totally dense subgroups. Thus
$\widehat{G}$ has also (at least) $m$-many totally dense subgroups by \cite[Lemma 4.1(c)]{comfsound82}.
Now apply again Proposition \ref{5.3.1}.
\epf
\mkp

\blem\label{cex1}
Let $\alpha\geq\cc$ be a cardinal, $B$ an Abelian group of bounded order with $|B|=\alpha$, and $H$  an infinite Abelian divisible group
with $2^{|H|}\leq \alpha$. Then the Abelian group $G:=B\times H$ with $|G|=\alpha$ admits at most $2^{2^{|H|}}<2^{2^\alpha}$ many SC-group
topologies.
\elem
\pf\
We have $\wG=\widehat{B}\times\wH$. By duality theory we get that $\widehat{B}$ is of bounded order, and $\wH$ is torsion-free. Hence $t\wG\cong\widehat{B}$
and $\wG/t\wG\cong\wH$. Now let $\tau$ be a \tb \ group topology such that every subgroup of $G$ is closed in $(G,\tau)$.
Then Lemma \ref{Lemma 1.} gives $t\wG\subseteq (G,\tau) \^ $. A result of Kakutani \cite{kaku43} implies $|\wH|=2^{|H|}\geq\cc$.
By Lemma \ref{D2} we get $r(\wH)=2^{|H|}$. Hence $\wH$ has exactly $2^{2^{|H|}}$-many subgroups. Thus there is the same number
of subgroups of $\wG$ containing $t\wH$.
\epf
\mkp

Now we are ready to give the announced negative answer to Question \ref{Q2}.

\bthm\label{cex2}
Let $\alpha$ and $\beta$ be infinite cardinals with $2^\beta\leq\alpha$. Then there is an Abelian group $G$ with $|G|=\alpha$ such that
$m=2^{2^\beta}<2^{2^{|G|}}$ holds for the number $m$ of SC-group topologies on $G$ which can be chosen in the following
ways:
\begin{enumerate}
\item [(a)] $G$ is a torsion group.
\item [(b)] $r_0(G)=\beta<\alpha$.
\end{enumerate}
\ethm

\pf\
We apply Lemma \ref{cex1}: Let $B$ as defined there and (a) $H:=\bigoplus_\beta Z(p^\infty)$ resp. (b) $H:=\bigoplus_\beta \QQ$.
Then, in both cases (a) and (b), the subgroup $H$ admits $2^{2^{\beta}}$-many SC-group topologies.
Indeed, it suffices to apply Lemma \ref{Prufer} in case (a), and Corollary \ref{newcor 1} in case (b), respectively.
Finally, applying Lemma \ref{nnewthm 2}, the proof is complete.
\epf
\mkp

Lemma \ref{nnewthm 2} is a main tool for the proof of

\bprp\label{nnewcor 1}
The following holds:
\begin{enumerate}
\item[(a)] Let $G$ be a countable Abelian group which is not reduced.
Then $G$ admits exactly $2^\cc$-many SC-group topologies.

\item[(b)] Every infinite countable divisible Abelian group admits exactly $2^\cc$-many SC-group topologies.

\item[(c)] Every Abelian group which is not reduced admits at least $2^\cc$-many SC-group topologies.

\item[(d)] Let $G$ be a divisible Abelian group with $\cf (|G|)>\omega$. Then $G$ admits exactly
$2^{2^{|G|}}$-many SC-group topologies.
\end{enumerate}
\eprp

\pf\
(a) $G$ contains $\mathbb{Q}$ or $Z(p^\infty)$ for some $p\in\PP$.
Then apply Lemma \ref{nnewthm 2}, Proposition \ref{khan2} and Theorem  \ref{newthm 1}.\\
(b) follows from (a).\\
(c) Apply Lemma \ref{nnewthm 2} and (b).\\
(d) By Corollary \ref{newcor 1} we may assume that $r_0(G)< |G|$. Then \cite[p. 85]{fuchsi} implies
$r(G)=\sum\limits_{p} r_p(G)$. Hence $|G|=\sum\limits_{p} r_p(G)$ by Lemma \ref{D2}.
Since $\cf (|G|)>\omega$, there is $p_0\in\PP$ such that $|G|=r_{p_0}(G)$. By \cite[p. 105]{fuchsi},
$|G|$ contains a subgroup isomorphic to $\bigoplus\limits_{|G|} \ZZ(p_0^\infty)$. Finally, it suffices to
apply Lemma \ref{Prufer} and Lemma \ref{nnewthm 2}.
\epf
\mkp

Concerning the problem whether a countable Abelian group which is not of bounded order admits exactly $2^\cc$-many SC-group topologies
we notice the following facts.\bkp

Set $K:=\prod_{n\in\NN}Z(p^n)$ and let $\mathbf p$ be a free ultrafilter on $\NN$ or, equivalently, a point in the
remainder $\beta\NN\setminus\NN$ of the Stone-\v Cech compactification $\beta\NN$ of $\NN$ furnished with the discrete topology. Every continuous function $f \colon\NN\to\TT$ can be extended to
a continuous function $\bar f\colon \beta\NN \to\TT$ (for more details see \cite[\S2]{Comfort1974}).
If $\sF$ is a filter on $\NN$ such that the filter with the basis $\{f(F):F\in\sF\}$ converges in $\TT$ to $t$,
we write $t=\lim\limits_\sF f(n)$. Now $\mathbf p$ converges in $\beta\NN$  to itself (\cite[section 6.5]{Gillman1976}).
Thus we have $$\bar f(\mathbf p)= \lim\limits_{\mathbf p} f(n).$$
For $n\in\NN$ we denote by $Z_{p^n}$ the copy of $Z(p^n)$ placed as a subgroup of $\mathbb{T}$ and
$$X_n\colon Z(p^n)\to \mathbb{T}$$ designates the canonical isomorphism of $Z(p^n)$ onto $Z_{p^n}$.\\
For all $x:=(x(n))\in K$ let $f_x \colon\NN\to\TT$ with $f_x(n):=X_n(x(n))$.
We define the function
$$\chi_\mathbf p\colon K\to \mathbb{T}$$ by $$\chi_\mathbf p[(x(n))]:=\lim\limits_{\mathbf p} f_x(n)$$ for all $(x(n))\in K$.
This means $\chi_\bp[(x(n))]=\overline{f_x}(\bp)$. Now $f_x(n)+f_y(n)=f_{x+y}(n)$ holds for all $n\in\NN$. Since $\NN$ is dense in $\beta\NN$, we get $\overline{f_x}(\bu)+\overline{f_y}(\bu)=\overline{f_{x+y}}(\bu)$ for all $\bu\in\beta\NN$. Hence $\chi_\bp$ is a character.

Let $U:=\oplus_{n\in\NN}Z(p^n)$. Then $\chi_\bp[U]=\{0\}$ holds. Since $U$ is dense in $K$, the character $\chi_\bp$ is discontinuous.

\blem\label{nlem0}{\rm(\cite[Proposition 5, p.~60]{Bourbaki1966})}
Let $\bu$ be an ultrafilter on a set $X$. If $A$ and $B$ are two subsets of $X$ such that $A\cup B\in\bu$, then either $A\in\bu$ or $B\in\bu$.
\elem
\mkp

\blem\label{nlem1} $\chi_\mathbf p[K]=\mathbb{T}$.
\elem
\pf\ Let $\alpha\in\mathbb{T}$ be arbitrary. Since $Z_{p^\infty}$ is dense in $\mathbb{T}$,
there is a sequence $(k_n/p^n)$ converging to $\alpha$ (here $k_n/p^n\in Z_{p^n}$ for all $n\in\NN$). Define
$x(n)=X_n^{-1}(k_n/p^n)$ for all $n\in\NN$. By using Lemma \ref{nlem0} we have
$$\chi_\mathbf p(x(n))=\lim\limits_{\mathbf p} X_n(x(n))= \lim\limits_{\mathbf p} k_n/p^n = \alpha.$$ Hence $\alpha\in\chi_p[K]$.
\epf
\mkp

The following assertion is well-known.

\blem\label{nlem2}
There is a copy $Q_p$ of $\mathbb{Q}$ such that $\mathbb{T}\cong Q_p\oplus L_p$ for some subgroup $L_p$ of $\mathbb{T}$.
\elem

\blem\label{nlem3}
Set $S_\mathbf p:=\chi_\mathbf p^{-1}[L_p]$.  Then $S_\mathbf p$ is a totally dense subgroup of $K$.
\elem
\pf\
By Lemma \ref{nlem1} and Lemma \ref{nlem2}  we have $K/S_\mathbf p\cong\mathbb{Q}$.
$\chi_\bp[U]=\{0\}$ implies that $S_\mathbf{p}$ is
dense in $K$. Thus it suffices to apply \cite[Corollary 3]{Khan1984}.\epf

\blem\label{nlem4}
$S_\mathbf p\not= S_\mathbf q$ if $\mathbf p\not= \mathbf q$.
\elem
\pf\ There is $A_\mathbf p\in \mathbf p$ with $A_\mathbf p\not\in \mathbf q$.
Hence $A_\mathbf q:=\NN\setminus A_\mathbf p\in\bq$ by Lemma \ref{nlem0}.
We take an element $(x(n))\in K$ satisfying the following requirements:\\
(1) $x(n)=0$ for all $n\in A_\mathbf q$;\\
(2) Take $\alpha\not= 0$ in $Q_p$.
Since $Z_{p^\infty}$ is dense in $\mathbb{T}$,
there is a sequence $(k_n/p^n)$ converging to $\alpha$ (here $k_n/p^n\in Z_{p^n}$ for all $n\in\NN$). Define
$x(n)=X_n^{-1}(k_n/p^n)$ for all $n\in\NN\setminus A_\mathbf q$.\\
We have $$\chi_\mathbf q[(x(n))]=\lim\limits_{\mathbf q} X_n(x(n))=0.$$
This implies that $(x(n))\in\chi_\mathbf q^{-1}(0)\subseteq S_\mathbf q$. On the other hand Lemma \ref{nlem0} implies
$$\chi_\mathbf p[((x(n))]=\lim\limits_{\mathbf p} X_n(x(n))=\lim\limits_{\mathbf p} k_n/p^n = \alpha\in Q_p.$$
Hence $(x(n))\notin \chi_p^{-1}(L_p)=S_\mathbf p$. This proves that $S_\mathbf p\not=S_\mathbf q$.
\epf

\bthm\label{nlem5}
$K$ contains $2^\cc$-many totally dense subgroups.
\ethm
\pf\
$\NN$ has $2^\cc$-many free ultrafilters by \cite[Corollary 7.4]{Comfort1974}.
Then apply Lemma \ref{nlem3} and Lemma \ref{nlem4}.
\epf
\mkp

Now Proposition \ref{5.3.1} yields

\bcor\label{ncor1}
$\bigoplus_{n\in\NN} Z(p^n)$ admits exactly $2^\cc$-many SC-group topologies.
\ecor
\mkp

We obtain

\bthm\label{nthm22}
Let $G$ be a torsion group such that $G$ has only finitely
many $p$-components $G_p$ and is not of bounded order.
Then $G$ admits at least $2^\cc$-many SC-group topologies.
\ethm
\pf \
There exists $p_0\in\PP$
with $|G_{p_0}|\geq\omega$ such that $L:=G_{p_0}$ is not of bounded order.
By \cite[(A.24)]{hewitt1963abstract} $L$ contains a subgroup $B$ such that
\begin{enumerate}
  \item[(a)] $B$ is isomorphic with a direct sum of cyclic $p_0$-groups;
  \item[(b)] $B$ is a pure subgroup of $L$;
  \item[(c)] the quotient group $L/B$ is divisible.
\end{enumerate}

We consider the following cases:
\begin{enumerate}
  \item[(1)] $B$ is not of bounded order. Then $B$ is isomorphic with a direct sum of cyclic $p_0$-groups, where infinitely many summands are pairwise different.
Hence it contains a subgroup $U$ being isomorphic with $H_{p_0}:=\bigoplus_{n\in\mathbb{N}}Z(p_0^n)$. Then Corollary \ref{ncor1} implies that $H_{p_0}$ has $2^\cc$-many SC-group topologies. By Lemma \ref{nnewthm 2} the same holds for $G$.
  \item[(2)] $B$ is of bounded order. Thus it is a proper subgroup of $L$. Since $B$ is a pure subgroup of $L$, then \cite[Theorem 27.5]{fuchsi}
implies that $B$ is a direct summand of $L$.
By (c) it follows that $L$ is not reduced. Thus Proposition \ref{nnewcor 1}(c) gives that $L$ admits $2^\cc$-many SC-group topologies.
Finally apply again Lemma \ref{nnewthm 2} to see that $G$ admits also $2^\cc$-many SC-group topologies.
\end{enumerate}
\epf
\mkp

{Now we can prove one of the main results in this section.}
\bkp

\noindent \textbf{Proof of Theorem \ref{Thm_C}}.
(a) Apply Theorem \ref{newthm 1}, Theorem \ref{nnnewthm 1}, Theorem \ref{nthm22} and Lemma \ref{nnewthm 2}.
(b) follows from (a).\epf

\section{Topologically simple groups}

We now consider the opposite direction from the section before, namely the existence of totally bounded group topologies without producing
closed non-trivial subgroups.
Topologically simple groups were introduced in Section 1. The following result gives a characterization of infinite Abelian
totally bounded groups which are topologically simple.

\bthm\label{ts1}
Let $(G,\tau)$ be an infinite Abelian totally bounded group
and let $S\leq \wG$ be its character group. Then the following assertions are equivalent:
\begin{enumerate}
  \item $(G,\tau)$ is topologically simple.
  \item $L\cap S=\{0\}$ for all proper closed subgroups $L$ of $\wG$. 
  \item All group topologies coarser than $\tau$ are Hausdorff or anti-discrete.
  \item Let $\mu$ be the totally bounded group topology on $S$ induced by the compact group topology of $\wG$. Then
  $(S,\mu)$ is topologically simple.
\end{enumerate}
\ethm
\pf\  \noindent (1)$\Rightarrow$ (2): Let $L$ be a proper closed subgroup of $\wG$. By Lemma \ref{closure}, we have $L=\AA(\wG, \AA(G,L))$, where $\AA(G,L)$ is not trivial. Since$ (G,\tau)$ is topologically simple, $\AA(G,L)$ is dense in it. By Theorem \ref{B} the intersection of $L$ and $S$ is trivial.

\noindent (2)$\Rightarrow$ (1):
Let $H$ be a non-trivial subgroup of $G$. Then $\AA(\wG,H)$ is a proper closed subgroup of $\wG$ by \cite{hewitt1963abstract} (23.24, Remarks (b) and (c)).
Hence (2) implies that the intersection of $\AA(\wG,H)$ and $S$ is trivial. By Theorem \ref{B}, $H$ is dense in $(G,\tau)$.

\noindent (1)$\Rightarrow$ (3) is trivial.

\noindent (3)$\Rightarrow$ (1). Let $N\subsetneq G$ be a closed subgroup of $(G,\tau)$. Since $N$ is closed,
the group $G/N$, equipped with the quotient topology $\tau_q$ defined by the canonical map $\pi\colon G\to G/N$, is
Hausdorff. Let $\mu:=\pi^{-1}(\tau_q)$ be the initial topology on $G$ defined by $\pi$. Since $N$ is $\mu$-closed, it follows
that $\mu$ is precompact, coarser than $\tau$ and not anti-discrete. By (3), this means that $N=\{0\}$.

\noindent (2)$\Rightarrow$ (4). Let $H$ be a closed subgroup of $(S,\mu)$ with $H\neq S$. Let $N$ be the closure of $H$ in $\wG$.
By $H\neq S$ and $H=N\cap S$ we have $N\neq \wG$. Now (2) implies $H=\{0\}$.

\noindent (4)$\Rightarrow$ (2). Let $L$ be a proper closed subgroup of $\wG$. Then $H:=L\cap S$ is a closed subgroup of $(S,\mu)$.
Assume $H=S$. Then $S\subseteq L$ and hence $L=\overline{L}=\wG$, which is impossible. Thus $H$ is a proper closed subgroup
of $(S,\mu)$. This implies $H=\{0\}$, since $(S,\mu)$ is topologically simple.
\epf

\bcor\label{ts2}
If $(G,\tau)$ is an infinite Abelian totally bounded topologically simple group, then:
\begin{enumerate}
	\item It is monothetic and torsion-free,
	\item its character group is monothetic and torsion-free,
	\item algebraically, it is a subgroup of the real numbers $\RR$,
	\item its weight is at most $\cc$.
\end{enumerate}
\ecor

\pf \ (1) follows directly from the definitions, and (2) from Theorem \ref{ts1}. \cite{hewitt1963abstract} (24.32) implies that $G$ is isomorphic
to a torsion-free subgroup of the (discrete) torus group $\TT$. Thus the torsion-free rank of $G$ is at most $\cc$. Considering its divisible hull
we get (3). By Theorem \ref{ts1}, the same holds for $S$. This implies the last assertion.\epf

\brem \label{tsintegers} Let $G$ be the group of integers. Then $\wG$ is the compact torus group $\TT$. By Corollary \ref{ts2} we have for every
dense subgroup $S$ of $\TT$: $(G,\tau_S)$ is topologically simple if and only if $S$ (as a subgroup of $\TT$) is topologically simple if and only if $S$ is a torsion-free subgroup of $\TT$.
Hence, there is a one-to-one correspondence between the totally bounded group topologies on $G$ which
produce topologically simple groups and the (dense) torsion-free
subgroups of $\TT$. Therefore, there are exactly $2^{\cc}$-many totally bounded group topologies on $G$ which produce topologically simple groups. The weight can be chosen to be $\cc$. (See Theorem \ref{Thm_D}
and Corollaries \ref{4.3} and \ref{7.4}, below.)
\erem

\brem \label{tsreals} Let $G$ be a subgroup of the real numbers. Then $G$ can be considered as a subgroup of the torus $\TT$.
Since $\TT$ has only finite proper closed subgroups,
the topology of $\TT $ induces a metrizable totally bounded group topology $ \tau $ on $ G $ such
that $ (G,\tau) $ is topologically simple.
\erem

\blem \label{tsinjective} Let $ G $ be an infinite Abelian group. $S$ is a subgroup of $\wG$ \st \ $G_S$ is  topologically simple  \sii \ every non-zero element of $S$ is injective.
\elem

\pf \ ($\Rightarrow$) If $g_1, g_2 \in G\setminus \{0\}, \phi \in S\setminus \{0\}$ with $g_1\neq g_2$, yet $\phi(g_1)=\phi(g_2)$, then $0\neq g_1-g_2 \in H:=\phi^{-1}[\{0\}] \implies H$ is a non-trivial proper subgroup of $G_S$, a contradiction.

($\Leftarrow$) Let $H$ be a closed subgroup of $G_S$, different than $G$. Then $G_S/H$ is \tb \ with character group $\AA(S,H)$ (Lemma \ref{2.1}). Notice that $\AA(S,H)\neq \{0\}$. If there were a non-zero element $h\in H$, then $\phi(h)=0$ for every non-zero element $\phi\in \AA(S,H)$, contradicting that all non-zero elements of $S$ are injective. It follows that $H= \{0\}$, as required.\epf
\mkp

\blem\label{injective}
Assume $1\leq\kappa\leq \cc$ and let $G:=\bigoplus_\kappa \ZZ$. Then there is a family $\sF$ of homomorphism $\phi :G \to\RR$ \st\ the following holds: $|\sF|=\cc$ and $\langle \sF\rangle\setminus\{0\}$ consists of injective functions.
\elem

\pf \ Let $\sB:=\{e_t: t\in\kappa\}$ be a generating independent subset of $G$, and $\sB':=\{u_t: t\in\cc\}$ a base of the vector space $\bigoplus_\cc\QQ$ over $\QQ$. Since $\cc=\sum_{i\in I}\kappa_i$ with $|I|=\cc$ and $\kappa_i=\kappa$ for all $i\in I$ by \cite[Proposition 4.4, Chapter III]{Levy1979}, $\RR$ is a union of
$\cc$-many pairwise disjoint subsets of size $\kappa$. Thus we get a family $\sE=\{f_i:{i\in I}\}$ with the following property $(\ast)$:
For all $i\in I$ the function $f_i:\kappa\to \cc$ is injective with $f_i[\kappa]\cap f_j[\kappa]=\emptyset$ for $i\neq j$.
It follows that $|\sE|=\cc$. For each $f\in\sE$ define $\phi_f:G\to\RR$ by $\phi_f(e_t):=u_{f(t)}$. Now let $\sF:=\{\phi_f :f\in\sE\}$. Then
$|\langle\sF\rangle|=\cc$. $\langle\sF\rangle\setminus\{0\}$ consists of injective functions: For, let $\psi=\sum_j a_j\phi_{f_j}\in\langle\sF\rangle\setminus \{0\}$ and $g=\sum_t b_te_t\in G\setminus \{0\}$.
Then $\psi(g)=(\sum_j a_j\phi_{f_j})(\sum_t b_te_t)=\sum_j\sum_t a_jb_tu_{f_j(t)}\neq 0$, since by property $(\ast)$ $\{u_{f_j(t)}\}$ is
linearly independent \wrt\ the vector space $\bigoplus_\cc\QQ$ over $\QQ$. Hence $\psi$ is injective.
\epf

\blem\label{extgt}
Let $G$ be a non-trivial subgroup of $\RR$. Then $G$ admits a \tb \ group topology $\tau$ such that $w(G,\tau)=\cc$ and $(G,\tau)$ is topologically simple.
\elem

\pf \ Let $m$ be the rank of $G$, and let $L:=\bigoplus_m\ZZ$. Since algebraically $\RR\subset\TT$, by Lemma \ref{injective} there is a family $\sF$ of homomorphisms
$\phi:L\to \TT$ such that $|\sF|=\cc$ and $\langle \sF\rangle\setminus\{0\}$ consists of injective functions. By \cite[(A.7)]{hewitt1963abstract} each $\phi$
can be extended to a character $\hat{\phi}:G\to \TT$. Let $\sM:=\{\hat{\phi}:\phi\in\sF\}$. Let $S:=\langle\sM\rangle$. Then $|\sM|=\cc$, and
$S$ is a subgroup of $\wG$ with $|S|=\cc$. Now each $\phi\in S \setminus\{0\}$ is injective: For, let $g\in G\setminus\{0\}$. Since $L$ is an essential subgroup of $G$, we have $\langle g\rangle\cap L\neq\{0\}$. Hence, there is $n \in \ZZ \setminus\{0\}$ with $ng\in L \setminus\{0\}$.
Since $\langle\sF\rangle \setminus\{0\}$ consists of injective functions, we get $0\neq \phi(ng)=n\phi(g)$. Hence $\phi(g)\neq0$.
Let $\tau$ be the precompact group topology on $G$ with character group $S$. By Lemma \ref{tsinjective} it is topologically simple, in particular totally bounded.
Since $|S|=\cc$, Lemma \ref{CO} implies $w(G,\tau)=\cc$.
\epf
\mkp

Now we are ready for the

\mkp
\noindent
\textbf{Proof of Theorem \ref{Thm_D}}.
Let $(G,\tau)$ be defined as in Lemma \ref{extgt}. Then $|S|=\cc$ holds for its character group $S$. By Theorem \ref{ts1}, $S$ is torsion-free. Lemma \ref{D2} implies $r(S)=\cc$. Then $\bigoplus_\cc\ZZ\subseteq S$. $\cc$ has $2^\cc$-many subsets of size $\cc$ by \cite[Proposition 5.2.14]{Ciesielski97}. Thus $S$ has $2^\cc$-many subgroups $U$ of cardinality $\cc$.
Let $\tau_U$ be the precompact group topology on $G$ with character group $U$. By Theorem \ref{ts1}(3),
it is totally bounded. Since $\tau_U\subseteq\tau$, $(G,\tau_U)$ is also topologically simple. $|U|=\cc$ implies $w(G,\tau_U)=\cc$ by Lemma \ref{CO}.

Let $\mu$ be any totally bounded group topology on $G$ such that $(G,\mu)$ is topologically simple. By Corollary \ref{ts2}(4), $w(G,\mu)\leq\cc$. Let $S'$
be the character
group of $(G,\mu)$. Then $S'\subseteq\wG$ with $|S'|\leq\cc$ and $|\wG|=2^{|G|}$. By \cite[Proposition 5.2.14]{Ciesielski97},
$\wG$ has exactly $(2^{|G|})^\cc=2^\cc$ many subsets of size $\leq\cc$. Hence $G$ has at most $2^\cc$-many totally bounded group topologies $\mu$
such that $(G,\mu)$ is topologically simple.
\epf
\mkp

Corollary \ref{7.4} below is an example of this Theorem.

\brem\label{Remark7.10}
The proof of Theorem \ref{Thm_D}
shows the following: Let $\sP$ be the power set of $\cc$, and let $G$ be a non-trivial subgroup of $\RR$.
Let $\sT\sS(G)$ be the set of all topologically simple \ totally bounded group topologies on $G$. Then there exists an order-preserving injection $f:\sP\to\sT\sS(G)$ such that $w(G,f[M])=|M|$ for all $M\in\sP$. There is $\tau\in\sT\sS(G)$ with $w(G,\tau)=\cc$, and there is $\sL\subseteq\sT\sS(G)$ with $|\sL|=2^\cc$ such that all $\mu\in\sL$ are coarser than $\tau$.

In \cite{Comfort1991a} the following notation is introduced: Let $\kappa$ and $\lambda$ be cardinals, and let $\sP(\kappa)$ be the power set of $\kappa$. Then $C(\kappa,\lambda)$ means that there is in $\sP(\kappa)$  a chain of length $\lambda$. By \cite{Comfort1991a} we have the following:
\begin{enumerate}
	\item {$C(\cc,\cc^+)$ holds by \cite[Corollary 1.7]{Comfort1991a},}
	\item $C(\cc,2^\cc)$ is not a theorem of ZFC {by \cite[Corollary 1.3]{Comfort1991a}.}
\end{enumerate}
Let $G$ be a non-trivial subgroup of $\RR$. Then:
\begin{enumerate}
	\item In $\sT\sS(G)$ there is a chain of length $\cc^+$,
	\item The following is not a theorem of ZFC : There is $\tau\in\sT\sS(G)$ with $w(G,\tau)=\cc$ and a chain $\sC\subseteq\sT\sS(G)$
with $|\sC|=2^\cc$ such that all $\mu\in\sT\sS(G)$ are coarser than $\tau$.
\end{enumerate}

\erem

\section{An Example: The integers}

{Let us consider the special case when $G=\ZZ$, the group of the integers. If $n \in \omega$, then $n\ZZ$ is a subgroup of $\ZZ$,
and conversely, if $H$ is a subgroup of $\ZZ$, there is a unique $n \in \omega$ \st \ $H=n\ZZ$.}

\blem\label{4.1} Let $S$ be a subgroup of $\wZ$. $k\ZZ$ is $\tau_S$-closed \sii \ $\frac{1}{k}\in S$.
\elem

\pf \ $(\Leftarrow)$ $(\frac{1}{k})^{-1}[\{0\}]=k\ZZ$.

\noindent
$(\Rightarrow)$ Having finite index in $\ZZ$, it follows that $ k\ZZ$ is $\tau_S$-open. Hence the map $\phi:\ZZ_S \to \ZZ_S/k\ZZ \simeq Z_k \to \< \frac{1}{k} \> \subset \wZ, 1 \mapsto \frac{1}{k}$ is $\tau_S$-continuous since it sends $k\ZZ$ to 0. Since $\phi = \frac{1}{k}$, it follows that $\frac{1}{k}\in S$. \epf

The following has been noticed in \cite[(3.5.4) and (3.5.5)]{Dik-Prod-Stoy1990}.

\bcor\label{4.2} Let $S$ be a subgroup of $\wZ$. $p^k\ZZ$ is $\tau_S$-closed for all $k\in \omega$ \sii \ $Z(p^\infty) \subseteq S$.
\ecor

\bcor\label{4.3} {Let $S$ be a subgroup of $\wZ$ \st \ $S \cap t\wZ=\{0\}$}. Then $\ZZ_S$ has no non-trivial closed subgroups.
\ecor

Since algebraically {$\TT = t\TT \oplus (\oplus_\cc \QQ)$} (\cite[(A.14)]{hewitt1963abstract}), we have the following two results.

\bcor\label{7.4} There exist $2^\cc$-many point-separating subgroups $S$ of $\wZ=\TT$ \st \ every subgroup of $\ZZ_S$ is dense.
\ecor

\pf \ If $X\subseteq \cc$, set $S_X:=\{0\} \oplus (\oplus_X \QQ)$. Obviously, $S_X$ is dense in $\TT$, hence it is point-separating. By the above corollary, $\ZZ_{S_X}$ has no non-trivial closed subgroups and {$X_1\neq X_2$ implies $S_{X_1} \neq S_{X_2}$}. \epf

The following result is a special case of Theorem \ref{nnnewthm 1}, but we give a short independent proof.

\bcor\label{7.5} There exist $2^\cc$-many point-separating subgroups $S$ of $\wZ=\TT$ \st \ every subgroup of $\ZZ_S$ is closed.
\ecor

\pf \ If $X\subseteq \cc$, set {$S_X:=t\TT \oplus (\oplus_X \QQ)$}. Obviously, $S_X$ is dense in $\TT$, hence it is point-separating.
By Corollaries \ref{A} or \ref{4.2}, every subgroup of $\ZZ_S$ is closed, and $X_1\neq X_2 \implies S_{X_1} \neq S_{X_2}$. \epf

The acronym lcm stands for {\em least common multiple}.

\blem\label{4.4a} If $n_1,n_2 \in \NN$, then $n_1\ZZ \cap n_2\ZZ=\mbox{\em lcm}(n_1,n_2)\ZZ$.
\elem

\pf \ Obviously $\supseteq$ holds. To see $\subseteq$, set $M=\mbox{lcm}(n_1,n_2)$ and assume that $n_1\ZZ \cap n_2\ZZ=t\ZZ$. Then there exist $z_1,z_2 \in \ZZ$ \st \ $n_1z_1=t=n_2z_2 \implies M|t$, as required. \epf

\blem\label{4.5} Let $F'=\{n_1,...,n_k\} \subset \NN$, $n \in \NN$ and $F=F'\cup \{n\}$. Then
\[\mbox{\em lcm }F=\mbox{\em lcm}~(\mbox{\em lcm }F',n). \]
\elem

\pf \ Set $M:=\mbox{lcm}(\mbox{lcm}F',n), N:=\mbox{lcm}F$, and $x:=\mbox{lcm}F'$. Then $M=\mbox{lcm}(x,n)$. Accordingly, $x|M$ and $n|M$. But $x=\mbox{lcm}F' \implies n_j|x\: (j=1,...,k)\implies n_j|M\: (j=1,...,k)$. Hence, $N|M$. On the other hand, $N=\mbox{lcm}F\implies n_j|N\: (j=1,...,k) \implies x|N$. Since $n|N$, it follows that $M|N$, as required. \epf

\blem\label{4.6} Let $C$ be a non-empty subset of $\NN$. Then
\[\left< \frac{1}{n}:n \in C\right > = \left< \frac{1}{\mbox{\em lcm }F}:F\in [C]^{<\infty} \right> , \]
as subgroups of $\wZ$, where $[C]^{<\infty}$ stands for all the finite subsets of $C$.
\elem

\pf \ $\subseteq$ Obvious (take $F=\{n\}$).

$\supseteq$ We prove that $\frac{1}{\mbox{lcm }F} \in$~LHS whenever $F:=\{n_1,...,n_k\}\subset C$. This obviously will prove the required contention. We use induction on $k$. Obviously true if $k=1$. Assume that $\frac{1}{\mbox{lcm }F} \in$~LHS, whenever $F$ has $t$ elements or less. Assume then that $F=\{n_1,...,n_t,n\}$. Set $F'=\{n_1,...,n_t\}$. As in the proof of Lemma \ref{4.5}, set $N:=\mbox{lcm}(\mbox{lcm}F',n)=\mbox{lcm}F$, and $x:=\mbox{lcm}F'$. By the induction hypothesis, $\frac{1}{x} \in$~LHS. In addition, let $d$ stand for the greatest common divisor 
of $x$ and $n$. Recall that there exists integers $a$ and $b$ \st \ $d=ax+bn$. Since $xn=Nd=N(ax+bn)$ we have that $\frac{1}{N}=\frac{ax+bn}{nx}=\frac{a}{n}+\frac{b}{x} \in$ LHS.
 \epf

\brem\label{4.4} Obviously, if $H_1$ and $H_2$ are closed subgroups of a \tg \ $G$, so is $H_1 \cap H_2$. When $G=\ZZ$, there are $n_1, n_2 \in \omega$ \st \ $H_1=n_1\ZZ$ and $H_2=n_2\ZZ$. It follows from Lemma \ref{4.4a} that $H_1 \cap H_2=$lcm$(n_1,n_2)\ZZ$. If $S$ is a subgroup of $\wZ$, set $C_S:=\{n\in \omega:n\ZZ \mbox{ is } \tau_S\mbox{-closed} \}$. By Lemma \ref{4.6} we have that $C_S$ is closed under taking lcms of finite subsets of $C_S$. Given $C \subseteq \NN$, set $S_C:= \< \frac{1}{n}: n \in C \> \subset \wZ$. If $\overline{C}:=\{n\in \NN:\frac{1}{n}\in S_C\}$, then $C \subseteq \overline{C}$, and since $S_C=S_{\overline{C}}$ (Lemma \ref{4.6}), we have that $\overline{C}$ is closed under taking lcms of its finite subsets. When $C \subseteq \NN$ is \st \ $C = \overline{C}$ we say that $C$ is {\em lcm-closed}. In Theorem \ref{Thm_B}
we saw the existence of a {greatest} subgroup $M S$ containing $S$ and such that $C_{MS}=C_S$. In the following result, we explicitly build $C_{MS}$ when $G=\ZZ$.
 \erem

\bthm\label{4.7} Let $S$ be a subgroup of $\wZ$. Then there exists a {greatest} dense subgroup $M{S}$ containing $S$ and \st \ $C_{MS}=C_S$.

Similarly, there exists a {smallest} subgroup $mS$ contained in $S$ \st \ $C_{mS}=C_S$. The group $mS$ is dense in $\wZ$ \sii \ $C_S$ is infinite. If not, there are no {smallest} {\em dense} subgroups $D \subseteq \wZ$ with $C_D=C_S$. Moreover, if $S$ is a torsion group,  then $mS=S$.
\ethm

\pf \ Set $mS:= \< \frac{1}{n}: n \in C_S \> \subset \wZ$. By Lemma \ref{4.1} $mS \subseteq S$, and by Remark \ref{4.4}, $C_{mS}=C_S$, hence $mS$ is as required. Obviously, $mS$ is finite \sii \ $C_S$ so is. Hence, $mS$ is dense in $\wZ$ \sii \ $C_S$ is infinite. Notice that $mS \subseteq$ $t\wZ$. If $\wZ \cong t\wZ \times F$ {(of course, algebraically, $F \cong \RR$)}, then $MS:=mS \times F$  is a maximal dense subgroup of $\wZ$ containing $S$ and \st \ $C_{MS}=C_S$.

Assume $C_S$ is finite and suppose $D$ is a dense subgroup of $\wZ$ \st \ $C_D=C_S$. Then $D/mS$ is infinite and contains an element $x+mS\: (x \in D)$ of infinite order. It follows that $\< x \> \cap mS = \{0\}$. Then $D':=\< mS \cup \{2x\} \> $ is a proper infinite subgroup of $D$ with $C_{D'}=C_S$.

Assume that $S$ is a torsion group and let $s \in S$. If $s\neq 0$, there are $m,n \in \NN$ relatively prime \st \ $s=\frac{m}{n}$. But then there is $k \in \NN$ \st \ $\frac{1}{n}=ks\in S \implies n \in C_S \implies \frac{1}{n} \in mS \implies s \in mS$. Hence, $mS=S$.
\epf
\bkp

Compare this result with Theorem \ref{Thm_B}.
\bkp

We can rewrite the above in terms of precompact group topologies:\bkp

\bcor\label{4.9d} Let $\tau$ be a precompact group topology on $\ZZ$. Then there exists a {greatest} \tb \ group topology $T\tau$ on $\ZZ$ producing precisely the same $\tau$-closed subgroups.

Similarly, there exists a {smallest}  precompact group topology ${t\tau}$ on $\ZZ$ producing precisely the same $\tau$-closed subgroups. Moreover, $t\tau$ is Hausdorff \sii \ there exist infinitely many $\tau$-closed subgroups. If not, then there are no {smallest} \tb \ group topologies $t$ with the same $\tau$-closed subgroups.
\ecor

\pf \ By the Comfort-Ross Theorem, there is a subgroup $S$ of $\wZ$ \st \ $\tau=\tau_S$. Applying Theorem \ref{4.7} to $S$, we obtain the result. \epf

\brem\label{4.10} Given $C \subseteq \NN$, set $S_C:= \< \frac{1}{n}: n \in C \> \subset \wZ$. Then we have $C_{S_C}=\overline{C}$ (see Remark \label (\ref{4.4}) and $mS_{\overline{C}}=S_{\overline{C}}=S_C$ (Lemma \ref{4.6}). The partially ordered set (under inclusion) $\sS_C:=[mS_C,MS_C]:=\{S \subseteq \wZ:S$ is a subgroup of $\wZ$ with $mS_C \subseteq S \subseteq MS_C \}$ {(see the proof of Theorem \ref{Thm_B})} 
has cardinality $2^\cc$ and maximal chains of length $\cc$. As mentioned before, $C$ is finite \sii \ $mS_C$ is finite.
Moreover, $mS_C$ is the only torsion group in (the poset) $[mS_C,MS_C]$. We also notice that $S_1, S_2 \in \sS_C$,
then $mS_C \subseteq S_1 \cap S_2$ \ie \ $\sS_C$ has the finite intersection property (fip), and $S_1, S_2 \in \sS_C \implies \< S_1 \cup S_2 \> \in \sS_C$.

We have thus assigned to any collection $\sC$ of subgroups of $\ZZ$ a subset $C$ of $\NN$. $\sC$ has fip \sii \ $C$ is lcm-closed. If this is the case, we have assigned a poset $\sS_\sC$ of subgroups of $\wZ$, \st \ if $S\in \sS_\sC$, then $\tau_S$ has precisely $\sC$ as its system of closed subgroups. The {smallest} element is a torsion-group which is infinite \sii \ $\sC$ is infinite. All of the elements in $\sS_\sC$, except the {smallest} one, have torsion-free elements.
\erem

\brem\label{4.11a} We can see the above in terms of precompact group topologies. Let $\tau$ be such a topology on $\ZZ$. By the Comfort-Ross Theorem, there is a subgroup $S$ of $\wZ$ \st \ $\tau=\tau_S$, and $\tau$ is Hausdorff \sii \ $S$ is infinite. Using the notation in Corollary \ref{4.9d}, the poset (under inclusion) $I_\tau:=[t\tau,T\tau]:=\{$precompact group topologies on $\ZZ$ between $t\tau$ and $T\tau \}$ has cardinality $2^\cc$ and maximal chains of length $\cc$. $t\tau$ is the only element in $I_\tau$ that could not be Hausdorff; it is \sii \ the system of $\tau$-closed subgroups is infinite. We also notice that if $\tau_1, \tau_2 \in I_\tau$, then $t\tau \subseteq \tau_1 \cap \tau_2$.
\erem

\brem\label{4.12} Let $S$ be a subgroup of $\wZ$. It should be obvious by now that it is its torsion part, $tS$, the one that produces the system of $\tau_S$-closed subgroups, \ie \ $k\ZZ$ is $\tau_S$-closed \sii \ $k\ZZ$ is $\tau_{tS}$-closed. Hence, setting $\tau=\tau_S$ in Remark \ref{4.11a}, there is a poset $I_\tau$ of precompact group topologies on $\ZZ$ producing the same system of $\tau_S$-closed subgroups. Again, $I_\tau$ has cardinality $2^\cc$, maximal chains of length $\cc$ and $\tau_{tS}$ is its {smallest} element which is Hausdorff \sii \ $tS$ is infinite. By now we can see that given subgroups $S_1, S_2$ of $\wZ$, they may or may not produce the same system of closed subgroups. They produce the same system of closed subgroups \sii \ $tS_1= tS_2$, \ie \ if they have the same torsion subgroups.
\erem

Given a precompact group topology $\tau$, and a subgroup of $\ZZ$, what is its  $\tau$-closure? The following results are to be compared with Theorem \ref{3.7}, Corollary \ref{5.2.1} and Theorem \ref{B}.

\blem\label{4.13} Let $S$ be a subgroup of $\wZ$ and let $k\ZZ$ be a subgroup of $\ZZ$. Then the $\tau_S$-closure of $k\ZZ$ equals $t\ZZ$ where $t\in C_S, t|k$, and $s\in C_S, s> t \implies s$ does not divide $k$.
\elem

\pf \ Assume that the $\tau_S$-closure of $k\ZZ$ equals $t\ZZ$. Then $t \in C_S$ and $k\ZZ \subseteq t\ZZ \implies t|k$. If $s\in C_S$, then $s\ZZ$ is $\tau_S$-closed, and if $s|k$, then $k\ZZ \subseteq s\ZZ \implies t\ZZ \subseteq s\ZZ \implies s|t \implies s\leq t$, as required. \epf

\bthm\label{4.14} Let $S$ be a subgroup of $\wZ$ and let $k\ZZ$ be a proper subgroup of $\ZZ$. Then $k\ZZ$ is  $\tau_S$-dense in $\ZZ$ \sii \ $k\not\in C_S$, and if $s\in C_S, s>1,$ then $s$ does not divide $k$.
\ethm

\pf \ {$(\Rightarrow)$} {Obviously $k\not\in C_S$. If there is $s\in C_S, s>1$ with $s|k$, then $k\ZZ \subseteq s\ZZ$ implies the $\tau_S$-closure of $k\ZZ$ is contained in $s\ZZ \neq \ZZ$, hence $k\ZZ$ wouldn't be $\tau_S$-dense in $\ZZ$.}

{$(\Leftarrow)$} Suppose $k\ZZ$ is not $\tau_S$-dense in $\ZZ$, and let $s\ZZ$ be its $\tau_S$-closure. Then $s\in C_S, s|k,$ and $ s>1$. If $s=k$, then $k\in C_S$, a contradiction. Then $s\neq k$, $k\not\in C_S$ and $s|k$, a contradiction. \epf

Some extreme cases:

\bcor\label{4.15} Let $S$ be a subgroup of $\wZ$. $S$ is a torsion-free subgroup of $\wZ$ \sii \ every proper subgroup of $\ZZ$ is $\tau_S$-dense in $\ZZ$.
\ecor

\bcor\label{4.16} Let $S$ be a subgroup of $\wZ$. $t\wZ \subseteq S$ \sii \ every proper subgroup of $\ZZ$ is $\tau_S$-closed in $\ZZ$.
\ecor

\section{Conclusion}
In this paper, we have dealt with the topological structure of totally bounded Abelian groups via their dual groups.
Every Abelian group $G$ has associated a {big} dual group $\wG$ that consists of the homomorphisms of $G$ into
the one-dimensional torus (the unit circle of the complex plane) in such a way that
each totally bounded Abelian group topology defined on $G$ coincides with the weak or initial topology defined by a \emph{separating}
subgroup of $\wG$.
Thus, each separating subgroup $S$ of $\wG$ defines a totally bounded group topology on $G$ and vice versa. We have exploited
this pairing: \emph{topology versus subgroup of $\wG$} in order to obtain real progress in the understanding of
totally bounded group topologies on Abelian groups. Among other results, we have calculated the number of totally bounded
group topologies that have a determined family of subgroups as closed subsets, and the number of totally bounded group topologies
that have all its subgroups closed.
We have also proved that the family of subgroups of $\wG$ that define the same collection of closed subgroups of $G$
always contain a {greatest} element {but in general not a smallest one.} Finally we have calculated the number of
topologically simple totally bounded group topologies for non-trivial subgroups of the real line, {the only possible groups having this property.}

\section{Acknowledgement}
The first author's research was partially supported by the Spanish Ministerio de Econom\'{\i}a y Competitividad, grant:
MTM/PID2019-106529GB-I00 (AEI/FEDER, EU) and by the Universitat Jaume I, grant UJI-B2019-
06.
\bigskip

\noindent {We thank the referee for her/his helpful comments.}

\bibliography{biblxclosd_subps}
\bibliographystyle{amsplain}


\end{document}